\documentclass[12pt, oneside]{amsart}
\usepackage[margin=1in]{geometry}
\usepackage{abstract}
\usepackage[bbgreekl]{mathbbol}
\usepackage{amssymb, amsmath, amsthm, mathrsfs, amscd, amsaddr}
\usepackage{subfig, color, graphicx}
\usepackage[dvipsnames]{xcolor}
\usepackage{microtype, siunitx, fancyhdr, polynom}
\usepackage{tikz, tikz-cd, physics, multicol, multirow, anyfontsize}
\usepackage{hyperref}
\usepackage{tabularx, makecell}

\hypersetup{colorlinks,linkcolor=blue,citecolor=cyan}

\usetikzlibrary{decorations.markings}

\newcommand{\m}[0]{\mathbb}

\newcommand{\B}[0]{\mathcal}
\newcommand{\G}[0]{\mathfrak}

\newcommand{\hb}[0]{\hbar}
\newcommand{\ov}[0]{\overline}
\newcommand{\wt}[0]{\widetilde}

\DeclareMathOperator{\Exp}{Exp}

\newcommand{\BoxedHom}{}
\def\BoxedHom[#1, #2](#3, #4)(#5){

    \draw[color=gray,thick] (#1 - #3/2, #2 - #4/2) -- (#1 + #3/2, #2 - #4/2) -- (#1 + #3/2, #2 + #4/2) -- (#1 - #3/2, #2 + #4/2) -- cycle;

    \draw (#1, #2) node  {#5};
}

\newcommand{\OverCross}{}
\def\OverCross[#1, #2](#3){

    \draw (#1, #2 - #3/4) -- (#1, #2) -- (#1 + #3, #2 + #3) -- (#1 + #3, #2 + 5*#3/4);
    \draw (#1 + #3, #2 - #3/4) -- (#1 + #3, #2) -- (#1 + #3/2 + #3/10, #2 + #3/2 - #3/10);
    \draw (#1 + #3/2 - #3/10, #2 + #3/2 + #3/10) -- (#1, #2 + #3) -- (#1, #2 + 5*#3/4);;
    
}

\newcommand{\UnderCross}{}
\def\UnderCross[#1, #2](#3){

    \draw (#1 + #3, #2 - #3/4) -- (#1 + #3, #2) -- (#1, #2 + #3) -- (#1, #2 + 5*#3/4);
    \draw (#1, #2 - #3/4) -- (#1, #2) -- (#1 + #3/2 - #3/10, #2 + #3/2 - #3/10);
    \draw (#1 + #3/2 + #3/10, #2 + #3/2 + #3/10) -- (#1 + #3, #2 + #3) -- (#1 + #3, #2 + 5*#3/4);
    
}

\newcommand{\HoroUnderCross}{}
\def\HoroUnderCross[#1, #2](#3){

    \draw (#1 + 5*#3/4, #2) -- (#1 + #3, #2) -- (#1, #2 + #3) -- (#1 - #3/4, #2 + #3);
    \draw (#1 - #3/4, #2) -- (#1, #2) -- (#1 + #3/2 - #3/10, #2 + #3/2 - #3/10);
    \draw (#1 + #3/2 + #3/10, #2 + #3/2 + #3/10) -- (#1 + #3, #2 + #3) -- (#1 + 5*#3/4, #2 + #3);
    
}

\newcommand{\HoroOverCross}{}
\def\HoroOverCross[#1, #2](#3){

    \draw (#1 + 5*#3/4, #2) -- (#1 + #3, #2) -- (#1 + #3/2 + #3/10, #2 + #3/2 - #3/10);
    \draw (#1 + #3/2 - #3/10, #2 + #3/2 + #3/10) -- (#1, #2 + #3) -- (#1 - #3/4, #2 + #3);
    \draw (#1 - #3/4, #2) -- (#1, #2) -- (#1 + #3, #2 + #3) -- (#1 + 5*#3/4, #2 + #3);
    
}

\newcommand{\TwoLabelHoro}{}
\def\TwoLabelHoro[#1, #2](#3)(#4, #5){

    \draw (#1, #2) -- (#1 + #3/2, #2);
    \draw (#1, #2 + #3) -- (#1 + #3/2, #2 + #3);

    \node at (#1 + #3/4, #2 + 5*#3/4) {\fontsize{\numexpr 12*#3}{\numexpr 12*12*#3/10}\selectfont #4};
    \node at (#1 + #3/4, #2 - #3/4) {\fontsize{\numexpr 12*#3}{\numexpr 12*12*#3/10}\selectfont #5};
}

\newcommand{\ThreeLabelHoro}{}
\def\ThreeLabelHoro[#1, #2](#3)(#4, #5, #6){

    \draw[color=gray,thick] (#1, #2) -- (#1 + #3/2, #2);
    \draw[color=gray,thick] (#1, #2 + #3) -- (#1 + #3/2, #2 + #3);
    \draw[color=gray,thick] (#1, #2 + 2*#3) -- (#1 + #3/2, #2 + 2*#3);
    
    \node at (#1 + #3/4, #2 + 9*#3/4) {#4};
    \node at (#1 + #3/4, #2 + 5*#3/4) {#5};
    \node at (#1 + #3/4, #2 - #3/4) {#6};
    
}

\allowdisplaybreaks

\tikzset{->-/.style={decoration={
  markings,
  mark=at position #1 with {\arrow{>}}},postaction={decorate}}}

\numberwithin{equation}{section}

\newtheorem{thm}[equation]{Theorem}

\newtheorem{lem}[equation]{Lemma}
\newtheorem{prop}[equation]{Proposition}

\newtheorem{conj}[equation]{Conjecture}

\author{Angus Gruen}

\address{Division of Physics, Mathematics and Astronomy, California Institute of Technology, 1200~E.~California Blvd., Pasadena, CA 91125, USA}

\email{agruen@caltech.edu}

\title{The $\G{sl}_N$ Symmetrically Large Coloured $R$ Matrix}

\begin{document}

\maketitle

\begin{abstract}
    For every knot $K$ and lie algebra $\G{g}$, there is a Gukov-Manolescu series denoted $F^{\G{g}}_K$ which serves as an analytic continuation of the quantum knot invariants associated to finite dimensional irreducible representations of $\G{g}$. There has been a great deal of work done on computing this invariant for $\G{g} = \G{sl}_2$ but comparatively less work has studied other lie algebras. In this paper we extend the large colour $R$ matrix from $\G{sl}_2$ to symmetrically coloured $\G{sl}_N$. This gives a definition for $F^{\G{sl}_N, sym}_K$ for positive braid knots and allows for predictions of $F^{\G{sl}_N, sym}_K$ for a much larger class of knots and links. It also provides further evidence towards a conjectural HOMFLY-PT analouge of $F_K$.
\end{abstract}

\tableofcontents

\section{Introduction}

Recently, there has been a great deal of interest in the knot invariant $F_K$ which acts as an analytic continuation of the coloured Jones polynomial. Recall the Melvin-Morton-Rozansky expansion \cite{MM, Roz, BNG}, for the coloured Jones polynomial which states that with $q = e^{\hb}$ and $x = q^n$:
\begin{align*}
    J_n(K; q = e^{\hb}) & = \sum_{m = 0}^{\infty}\sum_{j = 0}^{m} c_{K, m, j} n^j\hb^m
    \\ & = \frac{1}{\Delta_K(x)} + \sum_{k = 1}^{\infty} \frac{P_k(K; x)}{\Delta_K^{2k + 1}(x)} \frac{\hb^k}{k!}.
\end{align*}
Here, each $P_k(K; x)$ is a Laurent polynomial in $x$ and $\Delta_K(x)$ is the Alexander polynomial. It was conjectured in \cite{GM} that Borel resummation applied to the right hand side should give an integral two-variable series $F_K(x, q)$.

\begin{conj}[\cite{GM}, Conjecture 1.5]
    For every knot $K\subset S^3$, there exists an integer two-variable series $F_K(x, q)$ satisfying
    \[
        F_K(x, q) = \sum_{m = 0}^{\infty} f_m(q) x^m = \frac{1}{\Delta_K(x)} + \sum_{k = 1}^{\infty} \frac{P_k(x)}{\Delta_K^{2k + 1}(x)} \frac{\hb^k}{k!}. \label{eqn: sl2 MMR expansion} 
    \]
    Moreover, this series is annihilated by the quantum A-Polynomial \cite{Gar, Guk}
    \[
        \hat{A}_K(\hat{x}, \hat{y}, q)F_K(x, q) = 0  
    \]
\end{conj}

Note that, following \cite{EGGKPS} we express this conjecture using the \textit{positive reduced} expansion. In both the original paper, \cite{GM}, and elsewhere in the literature, $F_K$ is often given in the closely related \textit{balanced unreduced} expansion. To pass to this expansion\footnote{In the $\G{g} = \G{sl}_2$ case.} we should multiply by $x^{\frac{1}{2}} - x^{-\frac{1}{2}}$ and then perform the substitution $x^n \mapsto \frac{1}{2}(x^n - x^{-n})$.

Whilst in the initial paper, $F_K$ was only constructed for torus knots and $4_1$, over the course of the two papers \cite{Park2, Park3}, Park proved that $F_K$ could be defined for all links admitting a signed braid diagram\footnote{See \cite{Park3} Theorem 2 for a explicit definition of what a signed braid diagram is.}. In this work, Park also computed predictions for $F_K$ for a larger class of knots where the state sum method converged conditionally. So far these predictions also seem to be correct.

At the same time in \cite{Park1, EGGKPS} several extensions of $F_K$ were developed. As might be expected from its close relationship with the Jones polynomial, $F_K$ is naturally associated to the lie algebra $\G{sl}_2$. In \cite{Park1}, Park introduced the corresponding $F^{\G{g}}_K$ associated to other lie algebras $\G{g}$ and computed it for torus knots. This extension $F^{\G{g}}_K$ is more complicated as the corresponding coloured Jones polynomials $J^{\G{g}}_{\alpha}$ are now coloured not by an integer but a dominant integral weight. These weights live in a higher $n$ dimensional space and so $F^{\G{g}}_K$ is a series in $q$ and $x_1, \cdots, x_n$ as opposed to just $q$ and $x$.

Thus, instead of considering $F^{\G{g}}_K$ in full generality, in \cite{EGGKPS} a slight simplification denoted $F^{\G{sl}_N, sym}_K = F^N_K$ was considered which only concerns itself with symmetric representations of $\G{sl}_N$. Motivation for this comes from the HOMFLY-PT polynomial which interpolates the $\G{sl}_N$ Jones polynomials coloured by the fundamental representation. The idea was to develop a HOMFLY-PT analogue of $F_K$, which should be a three variable series $F_K(x, a, q)$ such that its specialisation $a = q^N$ recovers $F^N_K(x, q)$.

\begin{conj}[\cite{EGGKPS}, Conjecture 1]\label{conj: $a$-deformed $F_K$}
    For every knot $K\subset S^3$, there exists a~three-variable series $F_K(x,a,q)$ interpolating all $F^N_K$ in the following sense: 

    \begin{align*}
       & F_K(x,q^N,q)  = F^N_K(x,q),\\
       & \hat{A}_K(\hat{x},\hat{y},a,q)F_K(x,a,q) = 0.
    \end{align*}

    Moreover, it has the~following properties:
    \begin{align*}
    \label{Weyl_symmetry}
        F_K(x,a,q) & = F_K(a^{-1}x^{-1},a,q),\\
    F_K(x,1,q) & = \Delta_K(x), \\
     F_K(x,q,q) & = 1,\\
    \lim_{q\rightarrow 1}F_K(x,q^N,q) & = \frac{1}{\Delta_K(x)^{N-1}}.
    \end{align*}
    Its asymptotic expansion should agree with that of the~coloured HOMFLY-PT polynomials. That is, 
    \begin{equation*}
    \log F_K(e^{k\hbar},a,e^\hbar) = \log P_k(K;a,e^{\hbar})
    \end{equation*}
    as $\hbar$-series. 
\end{conj}

Closely related to this conjecture is the knots-quivers correspondence introduced in \cite{KRSS1, KRSS2} and adapted to $F_K$ in \cite{Kuch, EGGKPSS}. Roughly, this conjecture predicts that this $a$ deformed $F_K$ series, can be written as a quiver form
\begin{equation} \label{eq: Quiver Form}
    \exp(\frac{p(\log x,\log a)}{\hbar}) \sum_{\textbf{d}} \frac{(-q)^{\frac{1}{2}\textbf{d}M\textbf{d}^T}q^{\textbf{q}\cdot \textbf{d}}a^{\textbf{a}\cdot \textbf{d}}x^{\textbf{x}\cdot \textbf{d}}}{(q; q)_{\textbf{d}}}.
\end{equation}
Here $\textbf{x}, \textbf{a}$ are vectors of non negative integers, $M$ is a symmetric matrix of integers, $\textbf{q}$ should be a vector of half integers such that $\frac{1}{2}M_{ii} + \textbf{q}_i$ is always integral and $p$ is a polynomial of degree at most $2$. Hence, if $F^N_K(x, q)$ for some small values of $N$ is known, this equation form can be used as an ansatz to predict the full $a$ deformed $F_K$. Many of the results in this paper can be viewed providing evidence for both Conjecture \ref{conj: $a$-deformed $F_K$} and the knot-quiver correspondence as all the methods we use allow the construction of the aforementioned quiver forms for fixed specialisations $a = q^N$.

In this paper we prove rigorously the following theorem about the existence of $F^N_K(x, q)$ for a small but infinite class of knots.

\begin{thm}[$F_K$ for $\G{sl}_N$]\label{thm: FK for slN}
    Fix a positive (or negative) braid knot and a positive integer $N$. Then the knot invariant $F^N_K$ is well defined and satisfies

    \begin{align}
        \hat{A}_K(\hat{x},\hat{y},q^N,q)F^N_K(x,q) & = 0,
        \\ F_K(x^{-1}, a, q) & = F_K(q^N x, a, q)
        \\ F^N_K(x,q = e^{\hb}) & = \frac{1}{\Delta_K(x)^{N-1}} + \sum_{k = 1}^{\infty} \frac{R_k(x, N)}{\Delta^{N + 2k - 1}(x)} \hb^k. \label{eqn: MMR expansion} 
    \end{align}
\end{thm}

Additionally, whist currently unproven, we show how to extend the methods in \cite{Park3} to allow us to make predictions for $F^N_K(x, q)$ for a much larger class of knots. 

\begin{conj}[Inverse State Sum]\label{conj: FK for slN inv}
    The previous theorem extends to all homogeneous braid knots and more generally all links admitting a signed braid diagram (as in Theorem 2 of \cite{Park3}).
\end{conj}

Parts of this work, in particular the definition of the symmetrically large coloured $R$ matrix and a derivation for the $F_K$ series for the trefoil, were briefly mentioned in the authors collaboration \cite{EGGKPSS}. The main advances in this paper are a more rigorous treatment of this large colour $R$ matrix (In particular computing the prefactors correctly) as well as the proof of Theorem \ref{thm: FK for slN} and the extension of the methods in \cite{Park2, Park3} to the $\G{sl}_N$ case. Along side this, we present far more examples of explicit computations for $F^K_N$ for a collection of $K$ and $N$.

\subsection{Paper outline}
    In Section \ref{sec: Prerequisites} we cover some background representation theory for $\G{sl}_N$ and $U_q(\G{sl}_N)$. This is known to varying degrees by experts and can be skipped but is included for non experts and to fix notation and conventions.
    
    Section \ref{sec: large colour limit} explains how to extend \cite{Park2} for symmetrically coloured $\G{sl}_N$. It introduces the symmetrically coloured $\G{sl}_N$ $R$ matrix following the work of \cite{Bur} and shows how to take the symmetric large colour limit to pass to the $R$ matrix for a lowest weight Verma Module over $\m{C}(x^{\frac{1}{2}}, q^{\frac{1}{2}})$. It finishes by studying the classical limit of this representation and proving Theorem \ref{thm: FK for slN}.
    
    The remaining sections are devoted to computations and more experimental work. In particular section \ref{sec: Pos Braid Knots} computes $F^N_K(x, q)$ for a selection of positive braid knots with $N \in \{2, 3, 4\}$. It then explores predicting the $a$ deformations from these small $N$ computations. Additionally we introduce the stratified state sum which allows us to make predictions for $F^N_K(x, q)$ for some non positive braid knots.
    
    Finally, in Section \ref{sec: inverse state sum}, we introduce the inverted state sum method developed in \cite{Park2} and experimentally extend it to $\G{sl}_N$. This allows us to compute predictions for $F^N_K(x, q)$ for a much larger class of knots which in particular includes all knots admitting a signed braid diagram as in Theorem 2 from \cite{Park3}.

\section{Prerequisites} \label{sec: Prerequisites}

We start by briefly reviewing the structure of $\G{sl}_N$ and its symmetric representations. While this is well known, it will be useful to fix notation before we transition to the quantum group $U_q(\G{sl}_N)$.

\subsection{Symmetric representations of \texorpdfstring{$\G{sl}_N$}{sl(N)}}
    The Cartan matrix for the lie algebra $\G{sl}_N$ is the $(N - 1) \times (N - 1)$ square matrix
    \begin{align*}
        A_{\G{sl}_N} & = \begin{pmatrix}
            2 & -1 & 0 & 0 & \cdots\\
            -1 & 2 & -1 & 0 & \cdots \\
             &  & \ddots &  &  \\
            \cdots & 0 & -1 & 2 & -1\\
            \cdots & 0 & 0 & -1 & 2 \\
        \end{pmatrix}
    \end{align*}
    From this we can read off a generating set for $\G{sl}_N$ given by $(N - 1)$ $\G{sl}_2$ triples $(E_i, F_i, H_i)$ which satisfy the usual internal $\G{sl}_2$ relations as well as
    \begin{align*}
        [X_i, Y_j] & = 0 \quad \forall X, Y \in \{E, F, H\}, \ |i - j| > 1 \\
        [H_i, E_{i \pm 1}] & = -E_{i \pm 1} \\
        [H_i, F_{i \pm 1}] & = F_{i \pm 1} \\
        [E_i, F_{i \pm 1}] = [H_i, H_{i \pm 1}] & = 0 \\
        [E_i, [E_i, E_{i \pm 1}]] = [F_i, [F_i, F_{i \pm 1}]] & = 0.
    \end{align*}
    In particular from these relations, it can be verified that the map
    \begin{align*}
        E_i & \mapsto z_{i} \frac{\partial}{\partial z_{i + 1}} \quad \quad F_i \mapsto z_{i + 1} \frac{\partial}{\partial z_{i}} \\
        H_i & \mapsto z_{i} \frac{\partial}{\partial z_{i}} - z_{i + 1} \frac{\partial}{\partial z_{i + 1}}
    \end{align*}
    gives an action of $\G{sl}_N$ on $\m{C}[z_1, \cdots, z_N]$. As the action fixes the degree of monomials, the polynomial representation decomposes as
    \[
        \m{C}[z_1, \cdots, z_N] = \bigoplus_{k = 0}^{\infty} V_{N, k}
    \]
    where $V_{N, k}$ is the subrepresentation of homogeneous polynomials of degree $k$.
    \begin{lem}
        The representation $V_{N, k}$ is exactly the $k'$th symmetric representation of $\G{sl}_N$.
    \end{lem}
    
        

    We can define the action of $E, F, H$ triples corresponding to non simple positive roots using the adjoint action. Recall that the positive roots are indexed by pairs $1 \leq i \leq j \leq N - 1$ with $\alpha_{i, j} = \alpha_i + \alpha_{i + 1} + \cdots + \alpha_j$. Hence we define
    \begin{align*}
        E_{i, j} & = [E_i, [E_{i + 1}, [\cdots, [E_{j - 1}, E_j]]\cdots] = ad_{E_i}(ad_{E_{i + 1}}(\cdots ad_{E_{j - 1}}(E_j)) \cdots ) = z_i \frac{\partial}{\partial z_{j + 1}} \\
        F_{i, j} & = [F_j, [F_{j - 1}, [\cdots, [F_{i + 1}, F_i]]\cdots] = ad_{F_j}(ad_{F_{j - 1}}(\cdots ad_{F_{i + 1}}(F_i)) \cdots ) = z_{j+1} \frac{\partial}{\partial z_{i}}\\
        H_{i, j} & = H_i + H_{i + 1} + \cdots + H_j = z_i \frac{\partial}{\partial z_{i}} - z_{j + 1} \frac{\partial}{\partial z_{j + 1}}.
    \end{align*}
    In each case the final equality gives the induced action in the polynomial representation. As usual, we can equivalently think about these $V_{N, k}$'s as irreducible representations of the universal enveloping algebra $U(\G{sl}_N)$ and this is the starting point for describing the quantum analog of the above picture.

\subsection{Quantum analogs}
    While the following definitions are all elementary, there are a couple of differing conventions and so we use this section to explicitly fix the conventions for this paper. Define the quantum integers by
    \[
        [n]_q = \frac{q^{\frac{n}{2}} - q^{\frac{-n}{2}}}{q^{\frac{1}{2}} - q^{\frac{-1}{2}}}
    \]
    Using this, define the quantum factorial and derivative as:
    \[
        [n]_q! = \prod_{i = 1}^n [i]_q \quad \quad \left(\frac{\partial f}{\partial z}\right)_{q} = \frac{f(q^{\frac{1}{2}} z) - f(q^{-\frac{1}{2}} z)}{q^{\frac{1}{2}} z - q^{-\frac{1}{2}} z}.
    \]
	Clearly taking a limit as $q \to 1$ recovers the usual definitions and observe that, analogously to the classical case,
	\[
	    \left(\frac{\partial z^n}{\partial z}\right)_{q} = [n]_q z^{n - 1}.
	\]
    We also take a moment to introduce the $q$-Pochammer symbol depending on an integer $n$ or an integer vector \textbf{r}.
    \begin{align*}
        (x)_n = (x; q)_n & = \prod_{i = 0}^{n - 1} (1 - q^i x) \quad \quad (q)_{\textbf{r}} = \prod_{i} (q; q)_{r_i}.
    \end{align*}

\subsection{Symmetric representations of \texorpdfstring{$U_q(\G{sl}_N)$}{Uq(sl(N))}}
    We would like to construct a quantum analog $U_q(\G{sl}_N)$ for the universal enveloping algebra $U(\G{sl}_N)$. Following the Drinfeld Jimbo prescription, $U_q(\G{sl}_N)$ is generated by $N - 1$ tuples\footnote{Implicitly we are identifying $K_i$ with $K_{\omega_i}$ where $\omega_i$ is the fundamental weight dual to the $i$'th simple root.} $(E_i, F_i, K_i, K_i^{-1})$ which satisfy
    \begin{align*}
		K_0 & = 1 \\
		K_{\mu}K_{\lambda} & = K_{\mu + \lambda} \\
		K_{\lambda}E_iK_{\lambda}^{-1} & = q^{\frac{(\lambda, \alpha_i)}{2}} E_i \\
		K_{\lambda}F_iK_{\lambda}^{-1} & = q^{-\frac{(\lambda, \alpha_i)}{2}} F_i \\
		[E_i, F_j] & = \delta_{ij} \frac{K_i - K_i^{-1}}{q^{\frac{1}{2}} - q^{-\frac{1}{2}}} \\
		\sum_{n = 0}^{1 - A_{ij}} (-1)^n \frac{[1 - A_{ij}]_q!}{[1 - A_{ij} - n]_q![n]_q!} & E_i^n \ E_j \ E_i^{1 - A_{ij} - n} = 0 \\
		\sum_{n = 0}^{1 - A_{ij}} (-1)^n \frac{[1 - A_{ij}]_q!}{[1 - A_{ij} - n]_q![n]_q!} & F_i^n \ F_j \ F_i^{1 - A_{ij} - n} = 0
	\end{align*}
	We would like to quantize the polynomial representation described above. Nicely, this is almost immediate.
	\begin{prop} \label{prop: poly rep well defined}
	    The maps
	    \begin{align*}
            E_i & \mapsto z_{i} \left(\frac{\partial}{\partial z_{i + 1}}\right)_q \\
            F_i & \mapsto z_{i + 1} \left(\frac{\partial}{\partial z_{i}}\right)_q \\
            H_i & \mapsto z_{i} \frac{\partial}{\partial z_{i}} - z_{i + 1} \frac{\partial}{\partial z_{i + 1}} \\
            K_i & = q^{\frac{H_i}{2}}
        \end{align*}
        give a well defined representation of $U_q(\G{sl}_N)$ on $\m{C}(q^{\frac{1}{2}})[z_1, \cdots, z_N]$
	\end{prop}
	Verifying that this representation is well defined is easy but tedious so we will skip it here. It mostly boils down to the following pair of relations for quantum integers
	\begin{align*}
        [b + 1]_q[a]_q - [b]_q[a + 1]_q & = [a - b]_q \\
        [a + 2]_q - [2]_q[a + 1]_q + [a]_q & = 0
    \end{align*}
    Identically to the classical case, each generator fixes the overall degree of monomials and so this representation decomposes as a direct sum
    \[
        \m{C}(q^{\frac{1}{2}})[z_1, \cdots, z_N] = \bigoplus_{k = 0}^{\infty} V^q_{N, k}
    \]
    where $V^q_{N, k}$ is the irreducible subrepresentation of homogeneous polynomials of degree $k$. It will be handy to introduce a specific labelling for a basis $\{v_{\textbf{a}}\}$ of $V^q_{N, k}$. Here $\textbf{a}$ will denote an $n - 1$ tuple of integers $k \geq a_1 \geq \cdots \geq a_{n - 1} \geq 0$ and $v_{\textbf{a}}$ is related to the natural polynomial basis by
    \[
        v_{\textbf{a}} = \ket{a_0 = k, a_1, \cdots, a_{n - 1}, a_n = 0} = z_1^{a_0 - a_1} \cdots z_n^{a_{n - 1} - a_n}.
    \]
    The key feature of this basis is that the colour, $k$, only enters in a single place and so it will be easier to work with the $k \to \infty$ limit. With respect to this basis, our actions become
    \begin{align}
        E_i\cdot v_{\textbf{a}} & = [a_i - a_{i + 1}]_q v_{\textbf{a} - e_i} \nonumber \\
        F_i\cdot v_{\textbf{a}} & = [a_{i - 1} - a_i]_q v_{\textbf{a} + e_i} \label{eq: rep defn} \\
        K_i\cdot v_{\textbf{a}} & = q^{\frac{a_{i - 1} + a_{i + 1} - 2a_i}{2}}v_{\textbf{a}}. \nonumber
    \end{align}
    To find the actions for elements $E_{\alpha}, F_{\alpha}, K_{\alpha}$ corresponding to other positive roots we will need to replace the adjoint action with its quantised version \cite{Bur} given by
    \[
        ad^q_{X_{\alpha}}(X_{\beta}) =
        \begin{cases}
			q^{\frac{(\alpha, \beta)}{4}}X_{\alpha}X_{\beta} - q^{-\frac{(\alpha, \beta)}{4}}X_{\beta}X_{\alpha} & \text{if} \quad \alpha < \beta \\
			-ad^q_{X_{\beta}}(X_{\alpha}) & \text{if} \quad \alpha > \beta \\
			0 & \text{if} \quad \alpha = \beta.
		\end{cases}
    \]
    Here $X_{\alpha}, X_{\beta}$ are either both $E$'s or both $F$'s, the ordering\footnote{We use a slightly different but equivalent ordering to \cite{Bur}. Calling the ordering in \cite{Bur} $<_{B}$, the orderings are equivalent in the sense that $\alpha < \beta$ if and only if $\alpha <_{B} \beta$ or $(\alpha, \beta) = 0$.} is the reverse dictionary order $\alpha_{i, j} < \alpha_{i', j'}$ if $j < j'$ or $j = j'$ and $i < i'$ and the inner product is given by
    \[
		(\alpha_{i, j}, \alpha_{i', j'}) = \begin{cases}
			2 & (i, j) = (i', j')
			\\ 1 & i = i' \text{ or } j = j' \text{ but not both}
			\\ -1 & i = j' + 1 \text{ or } i' = j + 1
			\\ 0 & else.
		\end{cases}
	\]
	Note that if $(\alpha, \beta) = 0$ then $ad^q_{X_{\alpha}}(Y_{\beta})$ is also $0$. Then the action of the elements $E_{\alpha}, F_{\alpha}, K_{\alpha}$ are
    \begin{align*}
        E_{i, j} & = ad^q_{E_i}(ad^q_{E_{i + 1}}(\cdots ad^q_{E_{j - 1}}(E_j)) \cdots ) \\
        F_{i, j} & = ad^q_{F_j}(ad^q_{F_{j - 1}}(\cdots ad^q_{F_{i + 1}}(F_i)) \cdots )\\
        K_{i, j} & = K_{\alpha} = K_iK_{i + 1} \cdots K_j.
    \end{align*}
	Specialising to our representation, we find that there are some extra factors of $q$ on top of the obvious quantisation of the classical action
    \begin{align*}
		E_{i, j}\cdot v_{\textbf{a}} & = q^{\frac{j - i}{4} + \frac{a_j - a_i}{2}}[a_{j + 1} - a_j]_q v_{\textbf{a} - e_i - e_{i + 1} - \cdots - e_j}
		\\ F_{i, j}\cdot v_{\textbf{a}} & = q^{-\frac{j - i}{4} - \frac{a_j - a_i}{2}}[a_i - a_{i - 1}]_q v_{\textbf{a} + e_i + e_{i + 1} - \cdots + e_j}
		\\ K_{i, j}\cdot v_{\textbf{a}} & = q^{\frac{(a_j + a_i - a_{j + 1} - a_{i - 1})}{2}}v_{\textbf{a}}.
	\end{align*}
	
	\subsection{The quantum trace}
	    When we quantize the underlying algebra to $U_q(\G{sl}_N)$, this also quantizes the evaluation and co-evaluation\footnote{There are also $\overleftarrow{ev}^q_{N, k}$, $\overrightarrow{covev}^q_{N, k}$ maps but we ignore them here.} maps:
	    \begin{align}
	        \overrightarrow{ev}^q_{N, k}: & V^q_{N, k} \otimes (V^q_{N, k})^{*} \to \m{C}(q^{\frac{1}{2}}), \quad \quad v_{\textbf{i}} \otimes v_{\textbf{j}}^* \mapsto \left(\prod_{\alpha \in \Phi^+} K_{\alpha, \textbf{i}}\right) \delta_{\textbf{i}, \textbf{j}} \label{eq: evaluation}
	        \\ \overleftarrow{coev}^q_{N, k}: & \m{C}(q^{\frac{1}{2}})\to V^q_{N, k} \otimes (V^q_{N, k})^{*}, \quad \quad 1 \mapsto \sum_{\textbf{i}} v_{\textbf{i}} \otimes v_{\textbf{i}}^* \label{eq: co-evaluation}
	    \end{align}
	    Here by $K_{\alpha, \textbf{i}}$ we mean the eigenvalue satisfying $K_{\alpha} \cdot v_{\textbf{i}} = K_{\alpha, \textbf{i}} v_{\textbf{i}}$. On $V^q_{k, N}$ this factor is
	    \[
	        \left(\prod_{\alpha \in \Phi^+} K_{\alpha, \textbf{i}}\right) = q^{\frac{N - 1}{2}k - |\textbf{i}|}
	    \]
	    where $|\textbf{i}| = i_1 + \cdots + i_{N - 1}$ is the sum of the entries of \textbf{i}. Using these maps we can take the quantum trace of a function $f: V^q_{N, k} \to V^q_{N, k}$, denoted $\Tr^q_{N, k}(f) \in \m{C}(q^{\frac{1}{2}})$ by:
	    \[
	        \overrightarrow{ev}^q_{N, k} \circ (f \otimes 1) \circ \overleftarrow{coev}^q_{N, k}: \m{C}(q^{\frac{1}{2}}) \to \m{C}(q^{\frac{1}{2}}), \quad \quad 1 \mapsto \Tr^q_{N, k}(f).
	    \]
	    For an immediate application of this, let us compute the quantum dimension of $V^q_{N, k}$ given by the trace of the identity.
	    \[
	        \dim_q(V^q_{N, k}) = \Tr^q_{N, k}(\m{1}) = q^{\frac{N - 1}{2}k}\sum_{k \geq i_1 \geq \cdots \geq i_{N - 1} \geq 0} q^{-|\textbf{i}|} = q^{-\frac{N - 1}{2}k}\frac{(q^{k + 1} ;q)_{N - 1}}{(q;q)_{N - 1}}.
	    \]
	    Setting $a = q^N$ we recover (with a little manipulation) the HOMFLY-PT polynomial for the unknot from \cite{FGS}. If we additionally set $x = q^k$ we recover the fully unreduced $F_{0_1}(x, a, q)$ from \cite{EGGKPS}.
	    
	    This is exactly what we should expect as the unknot can be represented by a circle which is the graphical equivalent of trace of the identity. If we want to study more complicate knots, we need the ability to braid strands for which we introduce the $R_{\G{sl}_N}$ matrix. We also need to generalise this definition of a trace to functions on tensor products $f: (V^q_{N, k})^{\otimes i} \to (V^q_{N, k})^{\otimes i}$. This can by easily accomplished by composing $i$ copies of $\overrightarrow{ev}^q_{N, k}$ and $\overrightarrow{coev}^q_{N, k}$.
	
	\subsubsection{The $U_q(\G{sl}_N)$ $R$ Matrix} \label{sec: Generic R}
	
	    A general form for the $R$ matrix is given in \cite{Bur}:
        \begin{align}
            R_{\G{sl}_N} & = q^{\frac{f(\textbf{H})}{2}} \prod_{\alpha \in \Phi^+}^{<} \Exp_{q^{-1}}((1 - q^{-1}) K_{\alpha}^{\frac{1}{2}} E_{\alpha} \otimes K_{\alpha}^{-\frac{1}{2}} F_{\alpha}) \\
            \Exp_{q}(x) & = \sum_{r = 0}^{\infty} \frac{q^{\frac{r(r - 1)}{4}}x^r}{[r]_q!} \nonumber \\
            f(\textbf{H}) & = \sum_{i, j} a_{ij}^{-1} H_i \otimes H_j. \nonumber
        \end{align}
        Note that in comparing the above to \cite{Bur} we have slightly differing conventions which is what leads to the slightly unnatural definition for $\Exp_{q}$. To compute the product, we need to use the ordering for the roots which we defined a moment ago and looks like:
        \[
            \alpha_1 < \alpha_{1, 2} < \alpha_{2} < \alpha_{1, 3} < \alpha_{2, 3} < \alpha_{3} < \alpha_{1, 4} \cdots
        \]
        A nice feature of this choice of ordering is that we have a natural recursive formulation for the $R$ matrix:
        \[
           R_{\G{sl}_N} = q^{\frac{f_N(\textbf{H}) - f_{N-1}(\textbf{H})}{2}} R_{\G{sl}_{N - 1}} \times \prod_{i = 1}^{N - 1} \Exp_{q^{-1}}((1 - q^{-1}) K_{\alpha_{i, N}}^{\frac{1}{2}} E_{\alpha_{i, N}}K_{\alpha_{i, N}}^{-\frac{1}{2}} F_{\alpha_{i, N}}).
        \]
        Fixing a representation $V$, we can compute the matrix elements of $R_{\G{sl}_N}: V \otimes V \to V \otimes V$ and use this to produce the $V$-coloured quantum knot invariant.

\subsection{Quantum knot invariants}
    
    After the initial discovery of the Jones Polynomial in the mid 80's \cite{Jon}, it was soon fitted into a far more general framework (\cite{Wit} among others) which produced knot and link invariants from Quantum Groups and more generally solutions to the Yang-Baxter equation, \cite{Tur},
    \[
        \wt{R}_{23}\wt{R}_{12}\wt{R}_{23} = \wt{R}_{12}\wt{R}_{23}\wt{R}_{12}.
    \]
    The reason why solutions to the Yang-Baxter equation are important is that they naturally give rise to a representations of the braid group. Considering the representation category of $U_q(\G{sl}_N)$, we can define $\wt{R} = PR$ to be the $R$ matrix followed by the swap operator $P: v_{\textbf{i}} \otimes v_{\textbf{j}} = v_{\textbf{j}} \otimes v_{\textbf{i}}$. By $\wt{R}_{ij}$, we mean that $\wt{R}$ acts on the $i$'th and $j$'th components of the tensor product $V^{\otimes n}$.
    
    Thus once we fix a representation $V$ of $U_q(\G{sl}_N)$, $V^{\otimes n}$ carries a representation of the $n$-strand braid group $B_n$ given by $\sigma_i \mapsto \wt{R}_{i, i+1}$. Hence given a knot $K$ let $\beta_K$ be a braid whose right closure is the knot $K$. The simplest examples of this are the Unknot, Trefoil and Figure Eight as shown in Figure \ref{fig: Braid Closure 01, 31, 41}. Then, we can interpret the crossings as braid group elements and the cups and caps as  evaluation and co-evaluation maps. Thus each diagram becomes a map from the identity representation to itself which is simply a field element.
    
    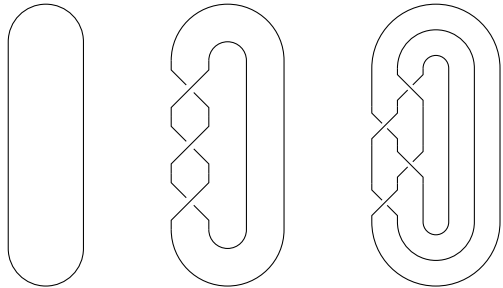
\begin{figure}[htp]
    	\centering
    	
    	\begin{tikzpicture}
    	    
    	    \def\size{2/4}
    	    
    	    
    	    \draw (0, 6*\size/4) arc(-180:0:3*\size/2) -- (3*\size, 24*\size/4) arc(0:180:3*\size/2);
    	    \draw (\size, 6*\size/4) arc(-180:0:\size/2) -- (2*\size, 24*\size/4) arc(0:180:\size/2);
    	    
    	    \OverCross[0, 7*\size/4](\size)
    	    
    	    \OverCross[0, 13*\size/4](\size)
    	    
    	    \OverCross[0, 19*\size/4](\size)
    	    
    	    
    	    \def\sizer{15/22 * \size}
    	    \def\xshft{16*\size/3}

    	    \draw (\xshft, 16*\sizer/4) -- (\xshft, 22*\sizer/4);
    	    \draw (\xshft + 2*\sizer, 22*\sizer/4) -- (\xshft + 2*\sizer, 28*\sizer/4);
    	    \draw (\xshft, 10*\sizer/4) arc(-180:0:5*\sizer/2) -- (\xshft + 5*\sizer, 34*\sizer/4) arc(0:180:5*\sizer/2) -- (\xshft, 28*\sizer/4);
    	    \draw (\xshft + \sizer, 10*\sizer/4) arc(-180:0:3*\sizer/2) -- (\xshft + 4*\sizer, 34*\sizer/4) arc(0:180:3*\sizer/2);
    	    \draw (\xshft + 2*\sizer, 16*\sizer/4) -- (\xshft + 2*\sizer, 10*\sizer/4) arc(-180:0:\sizer/2) -- (\xshft + 3*\sizer, 34*\sizer/4) arc(0:180:\sizer/2);
    	    
    	    \OverCross[\xshft, 11*\sizer/4](\sizer)
    	    
    	    \UnderCross[\xshft + \sizer, 17*\sizer/4](\sizer)
    	    
    	    \OverCross[\xshft, 23*\sizer/4](\sizer)
    	    
    	    \UnderCross[\xshft + \sizer, 29*\sizer/4](\sizer)
    	    
    	    
    	    \draw (\size - \xshft, \size) arc(-180:0:\size) -- (\size - \xshft + 2*\size, 13*\size/2) arc(0:180:\size) -- cycle;
    	    
        \end{tikzpicture}
        \captionsetup{singlelinecheck=off}
    	\caption[.]{Braid Closure diagrams for the Unknot $\beta_{0_1} = \m{1}$, Trefoil $\beta_{3_1} = \sigma_1^3$ and Figure Eight $\beta_{4_1} = \sigma_1\sigma_2^{-1}\sigma_1\sigma_2^{-1}$ knots.}
    	\label{fig: Braid Closure 01, 31, 41}
    \end{figure}
    
    Explicitly, choosing a representation $V$, the corresponding framed invariants will be:
    \[
        \ov{\B{P}}_V(\beta_{0_1}) = \Tr^q_V(\m{1}), \quad \ov{\B{P}}_V(\beta_{3_1}) = \Tr^q_{V^{\otimes 2}}(\wt{R}_{12}^3), \quad \ov{\B{P}}_V(\beta_{4_1}) = \Tr^q_{V^{\otimes 3}}(\wt{R}_{12}\wt{R}_{23}^{-1}\wt{R}_{12}\wt{R}_{23}^{-1}).
    \]
    \begin{lem}
        Fixing a finite dimensional representation $V$ of $U_q(\G{sl}_N)$ over $\m{C}(q)$, the map
        \[
            K \mapsto \ov{\B{P}}_V(\beta_{K})
        \]
        is a framed\footnote{Meaning a ribbon invariant, which obeys only the Reidemeister 2, 3 moves.} knot invariant.
    \end{lem}
    It remains to slightly adjust these definitions so as to remove the dependence on framing. First define the writhe of a braid as
    \[
        \omega(\beta) = |\beta_{+}| - |\beta_{-}|.
    \]
    Here $|\beta_{\pm}|$ denotes the number of positive/negative crossings in $\beta$. Then, in all cases we will consider\footnote{This occurs whenever the map corresponding to a Reidemeister 1 move is central meaning it has the form $f_V(q) \m{1}$.}, we can find a constant factor $f_V(q)$ such that for any $2$ braids $\beta_1, \beta_2$, whose closures represent the same knot,
    \[
        f_V(q)^{\omega(\beta_2)}\ov{\B{P}}_V(\beta_1) = f_V(q)^{\omega(\beta_1)} \ov{\B{P}}_V(\beta_2).
    \]
    Hence using this factor we get the full knot invariant:
    \[
        \wt{\B{P}}_V(K) = f_V(q)^{-\omega(\beta_K)}\ov{\B{P}}_V(\beta_K)
    \]
    where $\beta_K$ is any braid representative of $K$. Setting $V = V^q_{2, 2}$ this is exactly the celebrated unreduced Jones polynomial and for more general $N, k$ we get various unreduced coloured Jones polynomials. In practise, it is usually simpler to compute the reduced version of these invariants which, algebraically corresponds to dividing the unreduced invariant by the value of the invariant on the unknot. From the quantum group perspective, these reduced invariants come from the observation that, if we leave the left most strand open as in Figure \ref{fig: red Braid Closure 01, 31, 41}, braids become maps $V \to V$. If this map is central\footnote{Which it will be for the braids and representations we consider.}, it is given by $v \mapsto C_{\beta} v$ for some $C_{\beta}$ and so we can define the reduced trace $\wt{\Tr}^q_V(\beta) := C_{\beta}$. 
    
    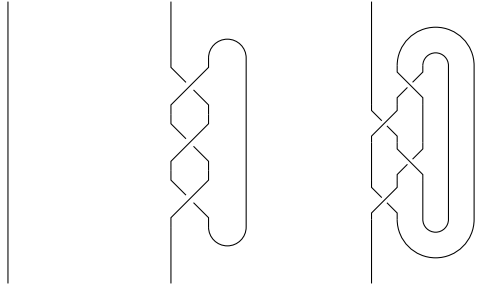
\begin{figure}[htp]
    	\centering
    	
    	\begin{tikzpicture}
    	    
    	    \def\size{2/4}
    	    
    	    
    	    \draw (0, 0) -- (0, 6*\size/4);
    	    \draw (0, 24*\size/4) -- (0, 30*\size/4);
    	    \draw (\size, 6*\size/4) arc(-180:0:\size/2) -- (2*\size, 24*\size/4) arc(0:180:\size/2);
    	    
    	    \OverCross[0, 7*\size/4](\size)
    	    
    	    \OverCross[0, 13*\size/4](\size)
    	    
    	    \OverCross[0, 19*\size/4](\size)
    	    
    	    
    	    \def\sizer{15/22 * \size}
    	    \def\xshft{16*\size/3}

    	    \draw (\xshft, 16*\sizer/4) -- (\xshft, 22*\sizer/4);
    	    \draw (\xshft + 2*\sizer, 22*\sizer/4) -- (\xshft + 2*\sizer, 28*\sizer/4);
    	    
    	    \draw (\xshft, 0) -- (\xshft, 10*\sizer/4);
    	    \draw (\xshft, 28*\sizer/4) -- (\xshft, 30*\size/4);

    	    \draw (\xshft + \sizer, 10*\sizer/4) arc(-180:0:3*\sizer/2) -- (\xshft + 4*\sizer, 34*\sizer/4) arc(0:180:3*\sizer/2);
    	    \draw (\xshft + 2*\sizer, 16*\sizer/4) -- (\xshft + 2*\sizer, 10*\sizer/4) arc(-180:0:\sizer/2) -- (\xshft + 3*\sizer, 34*\sizer/4) arc(0:180:\sizer/2);
    	    
    	    \OverCross[\xshft, 11*\sizer/4](\sizer)
    	    
    	    \UnderCross[\xshft + \sizer, 17*\sizer/4](\sizer)
    	    
    	    \OverCross[\xshft, 23*\sizer/4](\sizer)
    	    
    	    \UnderCross[\xshft + \sizer, 29*\sizer/4](\sizer)
    	    
    	    
    	    \draw (\size - \xshft, 0) -- (\size - \xshft, 30*\size/4);
    	    
        \end{tikzpicture}
        \captionsetup{singlelinecheck=off}
    	\caption[.]{Reduced Braid Closure diagrams for the Unknot $\beta_{0_1} = \m{1}$, Trefoil $\beta_{3_1} = \sigma_1^3$ and Figure Eight $\beta_{4_1} = \sigma_1\sigma_2^{-1}\sigma_1\sigma_2^{-1}$ knots.}
    	\label{fig: red Braid Closure 01, 31, 41}
    \end{figure}
    
    This finally leads us to the definition of the reduced invariant
    \[
        \B{P}_V(K) = f_V(q)^{-\omega(\beta_K)}\wt{\Tr^q_V}(\beta_K) = \frac{\wt{\B{P}}_V(K)}{\wt{\B{P}}_V(0_1)}.
    \]
    When $V = V^q_{N, k}$, the symmetric representations introduced above, this procedure yields a series of invariants known as the coloured HOMFLY-PT polynomials
    \[
        \B{P}_{V^q_{N, k}}(K) = P_k(K;a = q^N, q).
    \]
\section{The symmetric large colour limit} \label{sec: large colour limit}
    
    We follow the structure outlined in \cite{Park2} though we only focus on the lowest weight symmetric Verma module here. The picture for the highest weight module is analogous.

\subsection{An infinite Verma module for \texorpdfstring{$U_q(\G{sl}_N)$}{Uq(sl(N))}}

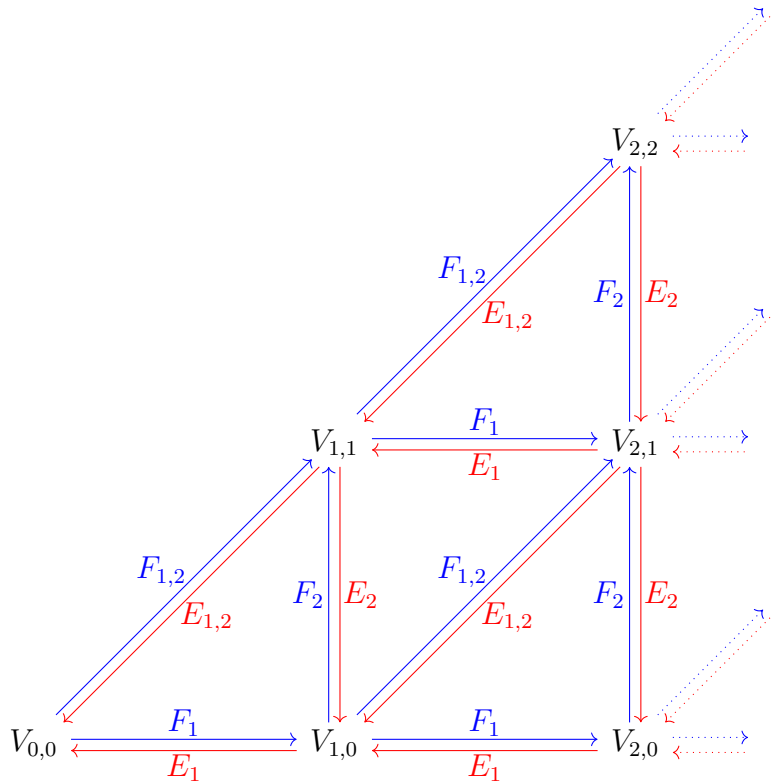
\begin{figure}[htp]
	\centering
	
	\begin{tikzpicture}
	    \tikzset{
	        font={\fontsize{12pt}{12}\selectfont}
	        }
	    \def\sc{4}
	    
	    \node[black] at (0, 0) {$V_{0, 0}$};
	    \node[black] at (\sc, 0) {$V_{1, 0}$};
	    \node[black] at (\sc, \sc) {$V_{1, 1}$};
	    \node[black] at (2*\sc, 0) {$V_{2, 0}$};
	    \node[black] at (2*\sc, \sc) {$V_{2, 1}$};
	    \node[black] at (2*\sc, 2*\sc) {$V_{2, 2}$};
	    
	    \draw[->, blue] (0.5, 0.075) -- (\sc - 0.5, 0.075);
	    \draw[<-, red]  (0.5, -0.075) -- (\sc - 0.5, -0.075);
	    \node[blue] at  (\sc/2, 0.3) {$F_1$};
	    \node[red] at   (\sc/2, -0.3) {$E_1$};
	    
	    \draw[->, blue] (\sc + 0.5, 0.075) -- (2*\sc - 0.5, 0.075);
	    \draw[<-, red]  (\sc + 0.5, -0.075) -- (2*\sc - 0.5, -0.075);
	    \node[blue] at  (1.5*\sc, 0.3) {$F_1$};
	    \node[red] at   (1.5*\sc, -0.3) {$E_1$};
	    
	    \draw[->, blue] (\sc + 0.5, \sc + 0.075) -- (2*\sc - 0.5, \sc + 0.075);
	    \draw[<-, red]  (\sc + 0.5, \sc - 0.075) -- (2*\sc - 0.5, \sc - 0.075);
	    \node[blue] at  (1.5*\sc, \sc + 0.3) {$F_1$};
	    \node[red] at   (1.5*\sc, \sc - 0.3) {$E_1$};
	    
	    \draw[->, blue] (\sc - 0.075, 0.3) -- (\sc - 0.075, \sc - 0.3);
	    \draw[<-, red]  (\sc + 0.075, 0.3) -- (\sc + 0.075, \sc - 0.3);
	    \node[blue] at  (\sc - 0.35, \sc/2) {$F_2$};
	    \node[red] at   (\sc + 0.35, \sc/2) {$E_2$};
	    
	    \draw[->, blue] (2*\sc - 0.075, 0.3) -- (2*\sc - 0.075, \sc - 0.3);
	    \draw[<-, red]  (2*\sc + 0.075, 0.3) -- (2*\sc + 0.075, \sc - 0.3);
	    \node[blue] at  (2*\sc - 0.35, 2) {$F_2$};
	    \node[red] at   (2*\sc + 0.35, 2) {$E_2$};
	    
	    \draw[->, blue] (2*\sc - 0.075, \sc + 0.3) -- (2*\sc - 0.075, 2*\sc - 0.3);
	    \draw[<-, red]  (2*\sc + 0.075, \sc + 0.3) -- (2*\sc + 0.075, 2*\sc - 0.3);
	    \node[blue] at  (2*\sc - 0.35, 1.5*\sc) {$F_2$};
	    \node[red] at   (2*\sc + 0.35, 1.5*\sc) {$E_2$};
	    
	    \draw[->, blue] (0.3, 0.4) -- (\sc - 0.3, \sc - 0.2);
	    \draw[<-, red]  (0.4, 0.3) -- (\sc - 0.2, \sc - 0.3);
	    \node[blue] at  (1.7, 2.3) {$F_{1,2}$};
	    \node[red] at   (2.3, 1.7) {$E_{1,2}$};
	    	    
	    \draw[->, blue] (\sc + 0.3, 0.4) -- (2*\sc - 0.3, \sc - 0.2);
	    \draw[<-, red]  (\sc + 0.4, 0.3) -- (2*\sc - 0.2, \sc - 0.3);
	    \node[blue] at  (1.5*\sc - 0.3, \sc/2 + 0.3) {$F_{1,2}$};
	    \node[red] at   (1.5*\sc + 0.3, \sc/2 - 0.3) {$E_{1,2}$};
	    
	    \draw[->, blue] (\sc + 0.3, \sc + 0.4) -- (2*\sc - 0.3, 2*\sc - 0.2);
	    \draw[<-, red]  (\sc + 0.4, \sc + 0.3) -- (2*\sc - 0.2, 2*\sc - 0.3);
	    \node[blue] at  (1.5*\sc - 0.3, 1.5*\sc + 0.3) {$F_{1,2}$};
	    \node[red] at   (1.5*\sc + 0.3, 1.5*\sc - 0.3) {$E_{1,2}$};
	    
	    \draw[dotted, ->, blue] (2*\sc + 0.5, 0.1) -- (2.5*\sc - 0.5, 0.1);
	    \draw[dotted, <-, red]  (2*\sc + 0.5, -0.1) -- (2.5*\sc - 0.5, -0.1);
	    
	    \draw[dotted, ->, blue] (2*\sc + 0.5, \sc + 0.1) -- (2.5*\sc - 0.5, \sc + 0.1);
	    \draw[dotted, <-, red]  (2*\sc + 0.5, \sc - 0.1) -- (2.5*\sc - 0.5, \sc - 0.1);
	    
	    \draw[dotted, ->, blue] (2*\sc + 0.5, 2*\sc + 0.1) -- (2.5*\sc - 0.5, 2*\sc + 0.1);
	    \draw[dotted, <-, red]  (2*\sc + 0.5, 2*\sc - 0.1) -- (2.5*\sc - 0.5, 2*\sc - 0.1);
	    
	    \draw[dotted, ->, blue] (2*\sc + 0.3, 0.4) -- (2.5*\sc - 0.3, \sc/2 - 0.2);
	    \draw[dotted, <-, red]  (2*\sc + 0.4, 0.3) -- (2.5*\sc - 0.2, \sc/2 - 0.3);
	    
	    \draw[dotted, ->, blue] (2*\sc + 0.3, \sc + 0.4) -- (2.5*\sc - 0.3, 1.5*\sc - 0.2);
	    \draw[dotted, <-, red]  (2*\sc + 0.4, \sc + 0.3) -- (2.5*\sc - 0.2, 1.5*\sc - 0.3);
	    
	    \draw[dotted, ->, blue] (2*\sc + 0.3, 2*\sc + 0.4) -- (2.5*\sc - 0.3, 2.5*\sc - 0.2);
	    \draw[dotted, <-, red]  (2*\sc + 0.4, 2*\sc + 0.3) -- (2.5*\sc - 0.2, 2.5*\sc - 0.3);
	    
    \end{tikzpicture}
    \captionsetup{singlelinecheck=off}
	\caption[.]{The bottom left corner of the lattice for symmetric representations of $U_q(\G{sl}_3)$. For the $k$'th symmetric representation, this lattice will stop at the line $V_{k, 0}, \cdots, V_{k, k}$ where for $i \leq j$, $V_{i, j}$ refers to the eigenspace spanned by:
	\[
	    v_{i, j} = \ket{k, i, j, 0} = z_1^{k - i}z_2^{i - j}.
	\]}
	\label{fig: Uqsl3 lattice}
    \end{figure}

    In the $U_q(\G{sl}_2)$ case, every irreducible representation can be thought of as a one dimensional lattice with nodes corresponding to an $H$ eigenspace of dimension $1$. Acting on an eigenspace by $E$ or $F$ has the effect of moving $1$ step either up or down this lattice. This picture extends to the irreducible symmetric representations of $U_q(\G{sl}_N)$ but, as might be expected, the lattice is no longer $1$ dimensional. For example when $N = 3$ we get the $2$ dimensional lattice illustrated in Figure \ref{fig: Uqsl3 lattice}. The triangular structure of this lattice comes simply from our choice of basis and note that for the $k$'th symmetric representation, the right hand boundary of this lattice are the eigenspaces $V_{k, i}$.
    
    The key algebraic feature of this choice of basis is that, the actions of $E_i, F_i, K_i$ given in Equation \eqref{eq: rep defn}, are mostly independant of the colour $k$. Indeed, $k$ appears only in the actions of $F_1, K_1$ as $q^{\pm \frac{k}{2}}$. Thus replacing $q^{k}$ by $x$ and extending the lattice in Figure \ref{fig: Uqsl3 lattice} to infinity we get a lowest weight Verma module $V^{x}_{N}$ of $\G{sl}_N$ over the field $\m{C}(q^{\frac{1}{2}}, x^{\frac{1}{2}})$.
    
    This module has basis
    \[
        v_{\textbf{a}} = \ket{a_1, \cdots, a_{n - 1}, a_n = 0}
    \]
    with $a_1 \geq \cdots \geq a_{n - 1} \geq 0$ with actions given by:
    \begin{align}
        E_i\cdot v_{\textbf{a}} & = [a_i - a_{i + 1}]_q v_{\textbf{a} - e_i} \nonumber \\
        F_1\cdot v_{\textbf{a}} = \frac{x^{\frac{1}{2}}q^{-\frac{a_1}{2}} - x^{-\frac{1}{2}}q^{\frac{a_1}{2}}}{q^{\frac{1}{2}} - q^{-\frac{1}{2}}} v_{\textbf{a} + e_1} \quad \quad F_i\cdot v_{\textbf{a}} & = [a_{i - 1} - a_i]_q v_{\textbf{a} + e_i} \nonumber \\
        K_1\cdot v_{\textbf{a}} = x^{\frac{1}{2}}q^{\frac{a_2 - 2a_1}{2}}v_{\textbf{a}} \quad \quad K_i\cdot v_{\textbf{a}} & = q^{\frac{a_{i - 1} + a_{i + 1} - 2a_i}{2}}v_{\textbf{a}} \nonumber
    \end{align}
    
    \begin{prop}
        The above definitions are well defined and give $V^{x}_{N}$ the structure of a infinite dimensional Verma Module.
    \end{prop}
    
    Most of the relations are identical to the ones required to check Proposition \ref{prop: poly rep well defined} and so follow from the identifies given below it. The only differences occur with relations involving $K_1$ and $F_1$. Lets explicitly check one of these:
    \[
        F_2^2F_1 - [2]_qF_2F_1F_2 + F_1F_2^2 = 0.
    \]
    Plugging in an arbitrary basis element we find:
    \begin{align*}
        \left(F_2^2F_1 - [2]_qF_2F_1F_2 + F_1F_2^2\right)v_{\textbf{a}} & = \frac{[a_1 - a_2]_q\left(x^{\frac{1}{2}}q^{-\frac{a_1}{2}} - x^{-\frac{1}{2}}q^{\frac{a_1}{2}}\right)}{(q^{\frac{1}{2}} - q^{-\frac{1}{2}})}
        \\ & \quad \quad \times \left([a_1 - a_2 + 1]_q - [2]_q[a_1 - a_2]_q + [a_1 - a_2 - 1]_q\right) v_{\textbf{a}}
        \\ & = 0
    \end{align*}
    The rest of the relations follow similarly.
    
    We can similarly extend the evaluation and evaluation maps to this module $V^{x}_{N}$. The co-evaluation map is identical to \eqref{eq: co-evaluation} but the evaluation becomes:
    \[
        \overrightarrow{ev}^{x}_{N}: V^{x}_{N} \otimes (V^{x}_{N})^{*} \to \m{C}(q^{\frac{1}{2}}, x^{\frac{1}{2}}), \quad \quad v_{\textbf{i}} \otimes v_{\textbf{j}}^* \mapsto x^{\frac{N - 1}{2}} q^{-|\textbf{i}|} \delta_{\textbf{i}, \textbf{j}}
    \]
    Again these combine to give a quantum trace and reduced quantum trace on $V^{x}_{N}$.
    
    This module $V_N^x$ should be thought of as a type of limit of the symmetric representations $V^q_{N, k}$. In particular, when we specialise $x = q^k$, $V^q_{N, k}$ sits inside $V_N^x$ as the irreducible subrepresentation containing $V_{\textbf{0}}$.
    
    It is worth noting that this property does not uniquely characterise $V^{x}_{N}$ and indeed there are other Verma modules we could have constructed. For example, there is an analogous generalisation of the highest weight Verma module described in \cite{Park2}. From the perspective of the polynomial representations we discussed earlier, this module corresponds to sending the exponent of $z_1$ to infinity while keeping the other exponents finite. This naturally leads to a host of other possible limits where, instead of $z_1$, the exponent of $z_j$ is sent to infinity. The highest weight module in \cite{Park2} and it's $\G{sl}_N$ extension corresponds to the $j = N$ case.

\subsection{The \texorpdfstring{$\G{sl}(N)$}{sl(N)} symmetrically large coloured R-matrix} \label{Sec: Large Colour R}
    
    In what follows by $R$ matrix we implicitly include the permutation operator and so technically these are all $\wt{R}$ matrices. Let us specialise the R matrix given in Section \ref{sec: Generic R} to the Verma module defined above. The result is an infinite summation over a collection of non negative integers $\textbf{r} = r_i^j$ with $1 \leq i \leq j \leq N - 1$. For $i \leq j \leq k \leq l$, define $\textbf{r}_{(i, j)}^{(k, l)}$ as the vector $(r_i^k, \cdots, r_i^l, r_{i + 1}^k, \cdots, r_j^l)$ and denote
	\[
		\textbf{r}^j  = \textbf{r}_{(1, j)}^{(j, j)} = (r_1^j, \cdots, r_j^j),
		\qquad\qquad \textbf{r}_j  = \textbf{r}_{(j, j)}^{(j, N - 1)} = (r_j^j, \cdots, r_j^{N - 1}).
	\]
	Letting $|\cdot |$ denote the $l^1$ norm we find that\footnote{Here we present the single coloured $R$ matrix which acts on $V^{x}_{N} \otimes V^{x}_{N}$. For link invariants, we should consider the multicoloured $R$ matrix which acts on $V^{x}_{N} \otimes V^{y}_{N}$. The matrix is essentially identical, we simply need to make the following replacement:
	\[
	    q^{\frac{(N - 1)\log_q(x)^2}{2N}}x^{-\frac{1}{2}(a_1 + b_1 + |\textbf{r}_1|)}(x q^{-b_1}; q^{-1})_{|\textbf{r}_1|} \mapsto q^{\frac{(N - 1)\log_q(x)\log_q(y)}{2N}}x^{-\frac{1}{4}(2b_1 + |\textbf{r}_1|)}y^{-\frac{1}{4}(2a_1 + |\textbf{r}_1|)}(y q^{-b_1}; q^{-1})_{|\textbf{r}_1|}.
	\]}:
	\begin{align} \label{eq: RMatrix slN}
		{}_{\G{sl}_N}R \ket{\textbf{a}, \textbf{b}} & = q^{\frac{(N - 1)\log_q(x)^2}{2N}} \sum_{\textbf{r} > 0} \frac{(-1)^{|\textbf{r}|}q^{C_N} x^{-\frac{1}{2}(a_1 + b_1 + |\textbf{r}_1|)}(x q^{-b_1}; q^{-1})_{|\textbf{r}_1|}(q^{a_1 - a_2 + |\textbf{r}_2|}; q^{-1})_{|\textbf{r}^1|}}{(q; q)_{\textbf{r}}}
		\\ & \quad \quad \times \prod_{j = 2}^{N - 1} (q^{b_{j - 1} - b_j}; q^{-1})_{|\textbf{r}_j|}(q^{a_j - a_{j + 1} + |\textbf{r}_{j + 1}|}; q^{-1})_{|\textbf{r}^j|} \ket{\textbf{a}', \textbf{b}'} \nonumber,
	\end{align}
	where
	\begin{align*}	
		C_N & = \frac{1}{2}\textbf{r} \cdot \textbf{r} + \textbf{a}\cdot M \cdot\textbf{b} +
		 \frac{1}{4} \sum_{j = 1}^{N - 1} |\textbf{r}^j|(a_{j + 1} + b_{j + 1} - 2) - \frac{1}{4} \left(\sum_{i = 2}^j r_i^j(a_{i - 1} + b_{i - 1})\right)
		\\ & \qquad + \sum_{i = 1}^j r_i^j\left(|\textbf{r}_{(i+1, j)}^{(j, N-1)}| + \frac{3}{4}(a_i - a_j) - \frac{1}{4}(b_i - b_j)\right),
		\\ M_{ij} & = \begin{cases}
			1 & i = j
			\\ -\frac{1}{2} & |i - j| = 1
			\\ 0 & \text{else.}
		\end{cases},
		\qquad a_i'  = b_i + |\textbf{r}_{(1, i)}^{(i, N - 1)}|.
		\qquad b_i'  = a_i - |\textbf{r}_{(1, i)}^{(i, N - 1)}|,
	\end{align*}
	Observe that the power of $x$ is always negative and, once we specialise to a particular vector $\textbf{r}$, the summand will always simplify to a polynomial. The second conclusion follows from the observation that, if we regroup some of the $q$-Pochammers, we find for each $j$
	\[
	    \frac{(q^{a_j - a_{j + 1} + |\textbf{r}_{j + 1}|}; q^{-1})_{|\textbf{r}^j|}}{(q; q)_{\textbf{r}^j}} = \frac{(q, q)_{a_j - a_{j + 1} + |\textbf{r}_{j + 1}| + |\textbf{r}^j|}}{(q; q)_{a_j - a_{j + 1} + |\textbf{r}_{j + 1}|} (q; q)_{\textbf{r}^j}}.
	\]
	The right hand side is simply a $q$-multinational coefficient and thus will be a polynomial in $q$. Also the prefactor $q^{\frac{(N - 1)\log_q(x)^2}{2N}}$, will be mostly cancelled out by the framing factor $f_N(x, q) = q^{\frac{(N - 1)\log_q(x)^2}{2N}} x^{\frac{N - 1}{2}}$. As we will always work in framing $0$, for computations we can replace this pre-factor by $x^{-\frac{N - 1}{2}}$ and ignore framing henceforth.
	
	To match the above description up with the R matrix given in \cite{Park2} define the matrix elements
    \begin{equation*}
    		{}_{\G{sl}_N}R_{\textbf{a}, \textbf{b}}^{\textbf{a}', \textbf{b}'} = \bra{\textbf{a}', \textbf{b}'} {}_{\G{sl}_N}R \ket{\textbf{a}, \textbf{b}}
    \end{equation*}
    where $\braket{\textbf{a}', \textbf{b}'}{\textbf{a}, \textbf{b}} = \delta_{\textbf{a}, \textbf{a}'}\delta_{\textbf{b}, \textbf{b}'}$. These matrix elements are $0$ unless $\textbf{a} + \textbf{b} = \textbf{b}' + \textbf{a}'$ in which case the summands will be non $0$ only when
    \[
        |\textbf{r}_{(1, i)}^{(i, N - 1)}| = a_i - b_i' = ab_i' - b_i.
    \]
    These conditions collapse the infinite summation to a finite sum and so each matrix element will be a polynomial in $x^{-1}, q$ and $q^{-1}$. We can similarly compute $R^{-1}$ matrix elements as:
    \[
        R^{-1} = P R|_{\substack{x \mapsto x^{-1} \\ q \mapsto q^{-1}}} P \quad \implies \quad {}_{\G{sl}_N}R^{-1}{}_{\textbf{a}, \textbf{b}}^{\textbf{a}', \textbf{b}'} = R{}_{\textbf{b}, \textbf{a}}^{\textbf{b}', \textbf{a}'}|_{\substack{x \mapsto x^{-1} \\ q \mapsto q^{-1}}}.
    \]
    Thus $R^{-1}$ matrix elements will be a polynomial in $x, q$ and $q^{-1}$.
    
    \subsection{The classical limit}
    
        Let's start by analysing the classical limit of this $R$ matrix. When we take $q \to 1$, the denominator has a $0$ of order $|\textbf{r}|$ and the numerator has a $0$ of order
        \[
            \sum_{j = 2}^{N - 1} |\textbf{r}^j| + \sum_{j = 1}^{N - 1} |\textbf{r}_j| = 2|\textbf{r}| - |\textbf{r}^1|.
        \]
        Hence in the $q \to 1$ limit, the only non $0$ terms occur when $|\textbf{r}| = |\textbf{r}^1|$ meaning $r_i^j = 0$ for $i \geq 2$. Thus we find (ignoring the prefactor for a moment): 
        \[
    		\lim_{q \to 1} {}_{\G{sl}_N}R \ket{\textbf{a}, \textbf{b}} = \sum_{\textbf{r}_1 > 0} (-1)^{|\textbf{r}_1|} x^{-\frac{1}{2}(a_1 + b_1 + |\textbf{r}_1|)}(1 - x)^{|\textbf{r}_1|}\prod_{j = 1}^{N - 1}\binom{a_j - a_{j + 1}}{r^j_1}\ket{\textbf{a}', \textbf{b}'}
    	\]
    	Passing to matrix elements, we find that the only non $0$ term in the summation is at $r^j_i = (a_i - a_{i + 1}) - (b'_i - b'_{i + 1})$ and so our $R$ matrix elements become:
    	\[
    		{}_{\G{sl}_N}R_{\textbf{a}, \textbf{b}}^{\textbf{a}', \textbf{b}'} = (-1)^{a_1 - b_1'} x^{-\frac{1}{2}(2a_1 + b_1 - b_1')}(1 - x)^{a_1 - b_1'}\prod_{j = 1}^{N - 1}\binom{a_j - a_{j + 1}}{b_i - b_{i+1}'}.
    	\]
    	This can be made simpler by considering a different labelling of our basis. Define $c_i = a_i - a_{i + 1}$ and $d_i =  b_i - b_{i + 1}$. Then with respect to this labelling (and reintroducing the prefactor modified by the framing) we have:
    	\[
    	    {}_{\G{sl}_N}R_{\textbf{c}, \textbf{d}}^{\textbf{c}', \textbf{d}'} = x^{\frac{N - 1}{2}}\prod_{i = 1}^{N - 1} (-1)^{c_i - d_i'} x^{-\frac{1}{2}(2c_i + d_i - d_i')}(1 - x)^{c_i - d_i'} \binom{c_i}{d_i'}.
    	\]
    	We immediately see that we have $n - 1$ non-interacting copies of the classical limit of ${}_{\G{sl}_2}R$. Hence, similarly to \cite{Park2}, if we compute the trace of a braid, after correcting for the framing we will recover $\frac{1}{\Delta_K(x)^{(N - 1)}}$. This proves the $\hbar = 0$ limit of property \eqref{eqn: MMR expansion} for knots where the $R$ matrix sum converges absolutely.
    	
    	As a brief side comment, note that in this $\hb = 0$ limit the theory is identical to the theory coming from $U_q(\G{sl}_2)^{N - 1}$. Similarly, if we study how the $R$ matrix acts on one particle states\footnote{These states are simplest to study in the \textbf{c}, \textbf{d} basis. In this basis, the one particle states are states where exactly one of the $c_i, d_i$ is $1$ and the rest are $0$.}, our theory is again identical to $U_q(\G{sl}_2)^{N - 1}$. The difference arises when $q$ is turned on where some multiparticle transitions occur only in the $U_q(\G{sl}_N)$ theory. This is one of the obstacles which currently prevents the generalization of the theorems in \cite{Park3} from $\G{sl}_2$ to $\G{sl}_N$.

    \subsection{Proof of Theorem \ref{thm: FK for slN}}
    
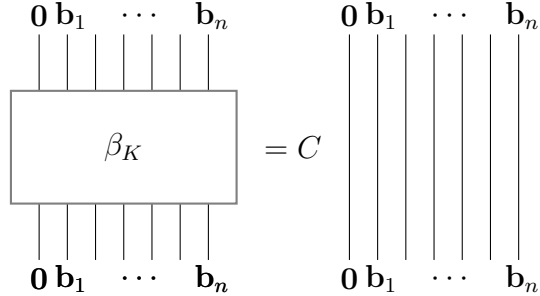
\begin{figure}[htp]
    	\centering
    	
    	\begin{tikzpicture}
    	    
    	    \def\size{3/4}
    	    \def\sft{\size/2}
    	    
    	    
    	    \node at (\sft, -0.25) {\textbf{0}};
    	    \node at (2*\sft + 0.05, -0.25) {$\textbf{b}_1$};
    	    \node at (9*\sft/2 + 0.05, -0.25) {$\cdots$};
    	    \node at (7*\sft + 0.05, -0.25) {$\textbf{b}_n$};
    	    
    	    \foreach \x in {1,...,7}{
    	        \draw (\x*\sft, 0)          -- (\x*\sft, \size);
    	        \draw (\x*\sft, 3*\size)    -- (\x*\sft, 4*\size);
    	        
    	        \draw (\x*\sft + 11*\sft, 0) -- (\x*\sft + 11*\sft, 4*\size);
    	    }
    	    
    	    \BoxedHom[2*\size, 2*\size](4*\size, 2*\size)($\beta_K$)
    	    
    	    \foreach \x in {0,1}{
    	        \foreach \y in {0,1}{
    	            \node at (\sft + \x*11*\sft,            -0.25 + \y/2 + \y*4*\size) {\textbf{0}};
    	            \node at (2*\sft + 0.05 + \x*11*\sft,   -0.25 + \y/2 + \y*4*\size) {$\textbf{b}_1$};
    	            \node at (9*\sft/2 + 0.05 + \x*11*\sft, -0.25 + \y/2 + \y*4*\size) {$\cdots$};
    	            \node at (7*\sft + 0.05 + \x*11*\sft,   -0.25 + \y/2 + \y*4*\size) {$\textbf{b}_n$};
    	        }
    	    }
 
    	    \node at (10*\sft, 2*\size) {$= C$};
    	    
        \end{tikzpicture}
        \captionsetup{singlelinecheck=off}
    	\caption[.]{The graphical picture defining the tensor elements $C = \bra{\textbf{0}, \textbf{b}_1, \cdots, \textbf{b}_n} \beta_K \ket{\textbf{0}, \textbf{b}_1, \cdots, \textbf{b}_n}$.}
    	\label{fig: beta matrix elements}
    \end{figure}

        We follow a similar outline as in \cite{Park2}. We first need to justify that the state sum converges absolutely for positive braid knots. To do this, let's look more closely at the $R$ matrix given in equation \eqref{eq: RMatrix slN}. In particular, observe that for each choice of $\textbf{r}$, the highest and lowest $x$ exponents which appear when we expand out the $q$-Pochammers will be
    	\begin{align}
    	    \text{Highest: } x^{-\frac{1}{2}(a_1 + b_1 - \textbf{r}_1)} & = x^{-\frac{1}{2}(b_1 + b_1')} \nonumber
    	    \\ \text{Lowest: }x^{-\frac{1}{2}(a_1 + b_1 + \textbf{r}_1)} & = x^{-\frac{1}{2}(a_1 + a_1')}. \nonumber
    	\end{align}
    	We see that if the incoming strands are labelled \textbf{a} and \textbf{b} then, regardless of the labelling on the outgoing strands, we will have an $x$ power of at most $x^{-\frac{b_1}{2}}$. As, when considering just the $x$ power, only the 1st component of the states \textbf{a}, \textbf{b}, \textbf{a}', \textbf{b}' appear, the situation is identical to the $\G{sl}_2$ case and so the argument in \cite{Park2} will easily carries across. That argument essentially proves the following lemma\footnote{Technically \cite{Park2} proves a slightly weaker lemma but the proof easily extends.}:
    	\begin{lem} \label{lem: min x degree}
    	    Let $\beta_K$ be an $n + 1$ strand braid representation of a knot $K$. Then
    	    \[
    	        \bra{\textbf{0}, \textbf{b}_1, \cdots, \textbf{b}_n} \beta_K \ket{\textbf{0}, \textbf{b}_1, \cdots, \textbf{b}_n}
    	    \]
    	    is a finite polynomial in $x^{-1}, q, q^{-1}$ with maximal $x$ coefficient
    	    \[
    	        x^{-\frac{1}{2}(b_{1, 1} + \cdots + b_{n, 1})}.
    	    \]
    	\end{lem}
    	For a graphical representation of $\bra{\textbf{0}, \textbf{b}_1, \cdots, \textbf{b}_n} \beta_K \ket{\textbf{0}, \textbf{b}_1, \cdots, \textbf{b}_n}$, see Figure \ref{fig: beta matrix elements}.
    	
    	The proof of this lemma follows from simple analysis of how weight can move around on a braid. Let us study the strand labelled by $\textbf{b}_i$ noting that, as we are dealing with a knot, $\sigma_i$ must appear at least once for all $i$. If $i = n$, then the first occurrence of $\sigma_n$ has bottom right strand $\textbf{b}_n$ and thus the corresponding $R$ matrix has maximal $x$ power $x^{-\frac{b_{n, 1}}{2}}$. If $i < n$ then there might be some number of $\sigma_{i + 1}$ before the first occurrence of $\sigma_i$.
    	
    	For any crossing going from state $\ket{\textbf{a}, \textbf{b}}$ to $\ket{\textbf{a}', \textbf{b}'}$ we immediately know $\textbf{a} \leq \textbf{a}' + \textbf{b}'$ and the corresponding $R$ matrix element has maximal $x$ power lower than $-\frac{b'_{1}}{2}$. Assume we encounter $j$ $\sigma_{i + 1}$ elements before we reach the first $\sigma_{i}$ element. Let the $k$'th $\sigma_{i + 1}$ go from $\ket{\textbf{a}_{k - 1}, \textbf{c}_{k}}$ to $\ket{\textbf{a}_k, \textbf{c}_{k}'}$ where $\textbf{a}_0 = \textbf{b}_i$. Then the overall maximum exponent power of $x$ coming from this chain of crossings is less than
    	\[
    	    -\frac{1}{2}(\textbf{c}_1' + \cdots + \textbf{c}_{j}' + \textbf{a}_j) \leq -\frac{1}{2}(\textbf{c}_1' + \cdots + \textbf{c}_{j - 1}' + \textbf{a}_{j - 1}) \leq \cdots \leq -\frac{\textbf{a}_0}{2} = -\frac{\textbf{b}_i}{2}
    	\]
    	Note that in this argument we have completely ignored the right incoming strand of $\sigma_{i + 1}$ and so we can freely apply this argument for all $i$ to conclude Lemma \ref{lem: min x degree}.
    	
    	Next, observe that if we fix the incoming labels, all summations over internal variables are finite. This is due to the fact that all labellings must be positive and, at every level, the sum of the labellings is $\textbf{b}_1 + \cdots + \textbf{b}_n$. This proves that the normalized reduced trace $f_{V^{x}_{N}}(q)^{-\omega(\beta_K)}\wt{\Tr}^q_{V^{x}_{N}}(\beta_K)$ converges absolutely for all positive braids $\beta_K$ to a series in $x^{-1}$ with coefficients, Laurent polynomials\footnote{From studying simple examples, we expect this can be improved to coefficients being simply polynomials in $q$ with $q^{-1}$ not appearing.} in $q$.
    	
    	Consider the properties of this series. In particular observe that, when we specialise $x = q^k$, we recover the quantum invariant $\B{P}_k(K, a = q^N, q)$. This follows from the fact that under the specialisation $x = q^k$ our module $V^x_N$ contains $V^q_{N, k}$ as the irreducible component containing $V_{\textbf{0}}$. As the open strand is coloured\footnote{A similar but more complicated argument still works if we had coloured the open strand by a different element.} $v_{\textbf{0}}$ and we are dealing with a knot\footnote{It should be possible to extend this argument to work for positive braid links as well.}, the trace restricts to the trace over this submodule $V^q_{N, k}$ which exactly computes the coloured Jones polynomial $\B{P}_{V_{N, k}}^q(K)$.
    	
    	As this is true for all $k$, it follows that $F^N_K(x, q) = f_N(x, q)^{-\omega(\beta_K)}\wt{\Tr}^q_{V^{x}_{N}}(\beta_K)$ is indeed an invariant of Positive Braid Knots which satisfies precisely the properties required by \ref{thm: FK for slN}.

\section{Our first examples}  \label{sec: Pos Braid Knots}
    
    We now turn to computing this invariant for some positive braid knots.

    \subsection{Positive braid knots}

    As can be seen in Table \ref{tab:pos braid knots}, there are $8$ positive braid knots with $10$ or less crossings. Of these, $6$ are torus knots, $4$ of which are $T(2, 2p + 1)$ torus knots. 
    
    \begin{table}
        \centering
        \begin{tabular}{ c|c|c|c }
        Knot & Braid & Torus Knot? \\ \hline
        $3_1^r$ & $\sigma_1^3$ & $T(2, 3)$ \\  
        $5_1^r$ & $\sigma_1^5$ & $T(2, 5)$ \\ 
        $7_1^r$ & $\sigma_1^7$ & $T(2, 7)$ \\  
        $8_{19}$ & $\sigma_1^3\sigma_2\sigma_1^3\sigma_2$ & $T(3, 4)$ \\  
        $9_1^r$ & $\sigma_1^9$ & $T(2, 9)$ \\  
        $10_{124}$ & $\sigma_1^5\sigma_2\sigma_1^3\sigma_2$ & $T(3, 5)$ \\  
        $10_{139}$ & $\sigma_1^4\sigma_2\sigma_1^3\sigma_2^2$ & No \\  
        $10_{152}$ & $\sigma_1^3\sigma_2^2\sigma_1^2\sigma_2^3$ & No \\  
        \end{tabular}
        \caption{Positive Braid knots with $10$ or less crossings. Data from \cite{LM}}
        \label{tab:pos braid knots}
    \end{table}
    
    Conjecture \ref{conj: $a$-deformed $F_K$} has been verified to a different extent for each of these knots. For $T(2, 2p + 1)$ torus knots, $F_K(x, a, q)$ is known in full generality \cite{EGGKPS}. For more general torus knots, $F_K(x, a, q)$ is not generally known (With the exception of $T(3, 4)$) but it is possible to compute $F_K(x, q^N, q)$ for any $N \in \m{N}$ via surgeries on plumbings \cite{Park1}. Finally, for non torus positive braid knots, only $F_K(x, q) = F_K(x, q^2, q)$ is only known \cite{Park2}.
    
    Thus we will first compute $F_K$ for a couple of torus knots as a cross check before focusing on $10_{139}$ and $10_{152}$.

    \subsubsection{Torus Knots}
    
    We start with $3$ torus knots, the Trefoil, Ciquefoil and $8_{19}$ knot. In each case, general $a$ deformed expressions are given in \cite{EGGKPS, EGGKPSS} which we can check against. As well as cross checking, this will give us an opportunity to briefly discuss how to efficiently perform these calculations. The main step is producing a labeled braid diagram\footnote{While usually these diagrams are drawn top to bottom, note we have drawn them left to right here.} as in Figure \ref{fig: Braid Diagrams 31, 51, 819}. To produce this diagram for a knot $K$, start with a braid whose closure is $K$. Give matching labels to corresponding left/right outgoing strands with the top stand labelled\footnote{In principal the top strand can carry any label as the choice of label will not affect the final answer but \textbf{0} is by the far the simplest for computations.} \textbf{0}.
    
    Next, at each crossing, ensure that the sum of incoming and outgoing labels are equal and that the label on the overstrand decreases. In particular this means that if the upper incoming/outgoing strand is labeled \textbf{0}, the diagonally opposite strand must also be labelled \textbf{0}. For simple knots such that the $3_1$ knot, this fixes all labels based off the external labels but in general this will introduce internal labellings which we need to sum over.
    
    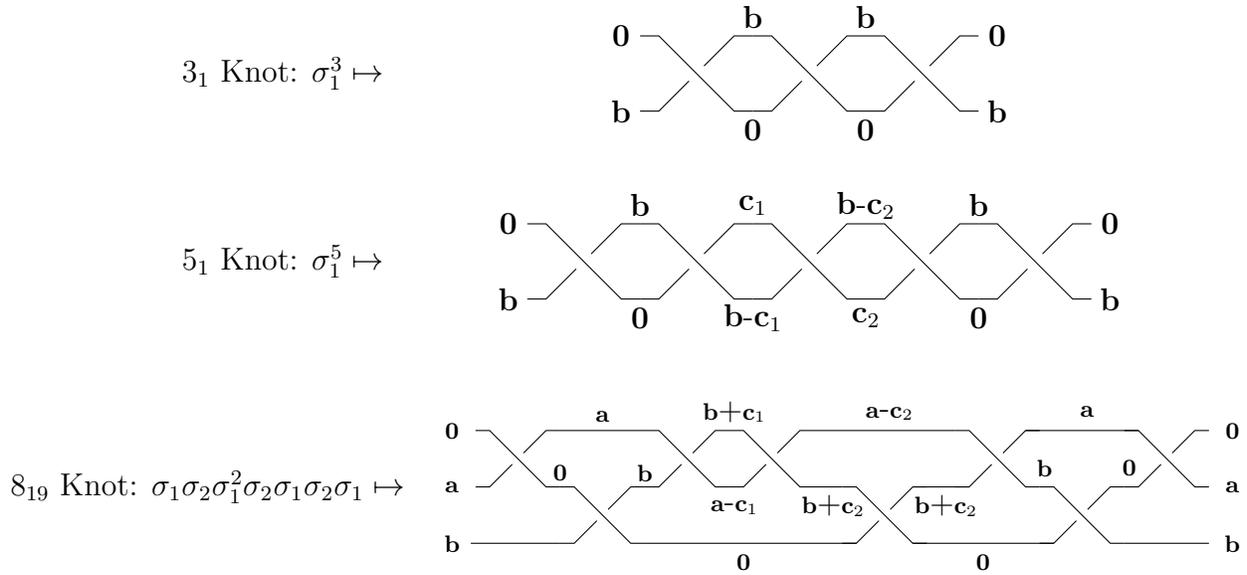
\begin{figure}[htp]
    	\centering
    	
    	\begin{tikzpicture}
    	    
    	    
    	    \node at (-5, 0.5) {$3_1$ Knot: $\sigma_1^3 \mapsto$};
    	
    	    \node at (-0.5, 1) {\textbf{0}};
    	    \node at (-0.5, 0) {\textbf{b}};
    	    
    	    \HoroUnderCross[0, 0](1)
    	    
    	    \node at (5/4, 1.25) {\textbf{b}};
    	    \node at (5/4, -0.25) {\textbf{0}};
    	    
    	    \HoroUnderCross[1.5, 0](1)
    	    
    	    \node at (11/4, 1.25) {\textbf{b}};
    	    \node at (11/4, -0.25) {\textbf{0}};
    	    
    	    \HoroUnderCross[3, 0](1)
    	    
    	    \node at (4.5, 1) {\textbf{0}};
    	    \node at (4.5, 0) {\textbf{b}};
    	    
    	    
    	    \def\yshft{2.5}
    	    \def\xshft{-1.5}
    	    
    	    \node at (-5, 0.5 - \yshft) {$5_1$ Knot: $\sigma_1^5 \mapsto$};
    	    
    	    \node at (\xshft - 1/2, 1 - \yshft) {\textbf{0}};
    	    \node at (\xshft - 1/2, 0 - \yshft) {\textbf{b}};
    	    
    	    \HoroUnderCross[\xshft, 0 - \yshft](1)
    	    
    	    \node at (5/4 + \xshft,  1.25 - \yshft) {\textbf{b}};
    	    \node at (5/4 + \xshft, -0.25 - \yshft) {\textbf{0}};
    	    
    	    \HoroUnderCross[\xshft + 1.5, 0 - \yshft](1)
    	    
    	    \node at (11/4 + \xshft,  1.25 - \yshft) {\textbf{c}\textsubscript{1}};
    	    \node at (11/4 + \xshft, -0.25 - \yshft) {\textbf{b}-\textbf{c}\textsubscript{1}};
    	    
    	    \HoroUnderCross[\xshft + 3, 0 - \yshft](1)
    	    
    	    \node at (17/4 + \xshft,  1.25 - \yshft) {\textbf{b}-\textbf{c}\textsubscript{2}};
    	    \node at (17/4 + \xshft, -0.25 - \yshft) {\textbf{c}\textsubscript{2}};
    	    
    	    \HoroUnderCross[\xshft + 4.5, 0 - \yshft](1)
    	    
    	    \node at (23/4 + \xshft,  1.25 - \yshft) {\textbf{b}};
    	    \node at (23/4 + \xshft, -0.25 - \yshft) {\textbf{0}};
    	    
    	    \HoroUnderCross[\xshft + 6, 0 - \yshft](1)
    	    
    	    \node at (\xshft + 7.5, 1 - \yshft) {\textbf{0}};
    	    \node at (\xshft + 7.5, 0 - \yshft) {\textbf{b}};
    	    
    	    
    	    \def\size{3/4}
    	    \def\xshftp{-2.25}
    	    \def\yshftp{}
    	    \def\fsizeone{9}
    	    \def\fsizetwo{10}
    	    
    	    \node at (-6, 0 - 2*\yshft) {$8_{19}$ Knot: $\sigma_1\sigma_2\sigma_1^2\sigma_2\sigma_1\sigma_2\sigma_1 \mapsto$};
    	    
    	    \node at (\xshftp - 0.5, \size - 2*\yshft) {\fontsize{\fsizeone}{\fsizetwo}\selectfont \textbf{0}};
    	    \node at (\xshftp - 0.5, 0 - 2*\yshft) {\fontsize{\fsizeone}{\fsizetwo}\selectfont \textbf{a}};
    	    \node at (\xshftp - 0.5, -\size - 2*\yshft) {\fontsize{\fsizeone}{\fsizetwo}\selectfont \textbf{b}};
    	    
    	    \draw (\xshftp - 1/4, -\size - 2*\yshft) -- (\xshftp + 5*\size/4, -\size - 2*\yshft);
    	    
    	    \HoroUnderCross[\xshftp, 0 - 2*\yshft](\size)
    	    
    	    \node at (\xshftp + 5*\size/4, \size/4 - 2*\yshft) {\fontsize{\fsizeone}{\fsizetwo}\selectfont \textbf{0}};
    	    
    	    \draw (\xshftp + 5*\size/4, \size - 2*\yshft) -- (\xshftp + 11*\size/4, \size - 2*\yshft);
    	    \node at (\xshftp + 2*\size, 5*\size/4 - 2*\yshft) {\fontsize{\fsizeone}{\fsizetwo}\selectfont \textbf{a}};
    	    
    	    \HoroUnderCross[\xshftp + 3*\size/2, -\size - 2*\yshft](\size)
    	    
    	    \node at (\xshftp + 11*\size/4, \size/4 - 2*\yshft) {\fontsize{\fsizeone}{\fsizetwo}\selectfont \textbf{b}};
    	    
    	    \draw (\xshftp + 11*\size/4, -\size - 2*\yshft) -- (\xshftp + 13*\size/2, -\size - 2*\yshft);
    	    \node at (\xshftp + 9*\size/2, -0.25 -\size - 2*\yshft) {\fontsize{\fsizeone}{\fsizetwo}\selectfont \textbf{0}};
    	    
    	    \HoroUnderCross[\xshftp + 6*\size/2, 0 - 2*\yshft](\size)
    	    
    	    \node at (\xshftp + 17*\size/4, 1/4 + \size - 2*\yshft) {\fontsize{\fsizeone}{\fsizetwo} \textbf{b}+\textbf{c}\textsubscript{1}};
    	    \node at (\xshftp + 17*\size/4, -1/4 - 2*\yshft) {\fontsize{\fsizeone}{\fsizetwo} \textbf{a}-\textbf{c}\textsubscript{1}};
    	    
    	    \HoroUnderCross[\xshftp + 9*\size/2, 0 - 2*\yshft](\size)
    	    
    	    \draw (\xshftp + 23*\size/4, 0 - 2*\yshft) -- (\xshftp + 25*\size/4, 0 - 2*\yshft);
    	    \node at (\xshftp + 24*\size/4, -0.25 - 2*\yshft) {\fontsize{\fsizeone}{\fsizetwo} \textbf{b}+\textbf{c}\textsubscript{2}};
    	    
    	    \draw (\xshftp + 23*\size/4, \size - 2*\yshft) -- (\xshftp + 33*\size/4, \size - 2*\yshft);
    	    \node at (\xshftp + 28*\size/4, 0.25 + \size - 2*\yshft) {\fontsize{\fsizeone}{\fsizetwo} \textbf{a}-\textbf{c}\textsubscript{2}};
    	    
    	    \HoroUnderCross[\xshftp + 13*\size/2, -\size - 2*\yshft](\size)
    	    
    	    \draw (\xshftp + 15*\size/2, -\size - 2*\yshft) -- (\xshftp + 20*\size/2, -\size - 2*\yshft);
    	    \node at (\xshftp + 35*\size/4, -0.25 -\size - 2*\yshft) {\fontsize{\fsizeone}{\fsizetwo}\selectfont \textbf{0}};
    	    
    	    \draw (\xshftp + 15*\size/2, 0 - 2*\yshft) -- (\xshftp + 17*\size/2, 0 - 2*\yshft);
    	    \node at (\xshftp + 16*\size/2, -0.25 - 2*\yshft) {\fontsize{9}{10} \textbf{b}+\textbf{c}\textsubscript{2}};
    	    
    	    \HoroUnderCross[\xshftp + 17*\size/2, 0 - 2*\yshft](\size)
    	    
    	    \node at (\xshftp + 39*\size/4, 0.25 - 2*\yshft) {\fontsize{\fsizeone}{\fsizetwo} \textbf{b}};
    	    
    	    \draw (\xshftp + 19*\size/2, \size - 2*\yshft) -- (\xshftp + 23*\size/2, \size - 2*\yshft);
    	    \node at (\xshftp + 21*\size/2, \size + 0.25 - 2*\yshft) {\fontsize{\fsizeone}{\fsizetwo} \textbf{a}};
    	    
    	    \HoroUnderCross[\xshftp + 20*\size/2, -\size - 2*\yshft](\size)
    	    
    	    \node at (\xshftp + 45*\size/4, 0.25 - 2*\yshft) {\fontsize{\fsizeone}{\fsizetwo} \textbf{0}};
    	    
    	    \draw (\xshftp + 45*\size/4, -\size - 2*\yshft) -- (\xshftp + 51*\size/4, -\size - 2*\yshft);
    	    
    	    \HoroUnderCross[\xshftp + 23*\size/2, 0 - 2*\yshft](\size)

    	    \node at (\xshftp + 25*\size/2 + 0.5,  \size - 2*\yshft)     {\fontsize{9}{10}\selectfont \textbf{0}};
    	    \node at (\xshftp + 25*\size/2 + 0.5,      0 - 2*\yshft)     {\fontsize{9}{10}\selectfont \textbf{a}};
    	    \node at (\xshftp + 25*\size/2 + 0.5, -\size - 2*\yshft)     {\fontsize{9}{10}\selectfont \textbf{b}};
    	    
        \end{tikzpicture}
        \captionsetup{singlelinecheck=off}
    	\caption[.]{Labeled Braid diagrams for the Trefoil ($3_1$) and Cinquefoil ($5_1$) and $8_{19}$ knots.}
    	\label{fig: Braid Diagrams 31, 51, 819}
    \end{figure}
    
    Finally, we apply the quantum trace, reading the $R$ matrix elements directly off the diagram to get:
    \begin{equation}
        \sum_{\textbf{a}_1, \cdots, \textbf{a}_n, \textbf{c}_1, \cdots \textbf{c}_m} \prod_{i = 1}^n x^{\frac{N - 1}{2}} q^{-|\textbf{a}_i|} \prod_{\alpha \in \text{crossings}} {{}_{\G{sl}_N}R}_{\textbf{i}_\alpha, \textbf{j}_\alpha}^{\textbf{i}_\alpha', \textbf{j}_\alpha'}. \label{eq: r mat braid sum}
    \end{equation}
    This will give a series in $x^{-1}$ so as a final step we apply Weyl symmetry, sending $x^{-1} \to q^Nx$ to get a series in $x$.
    
    Let's start by applying this to the right handed trefoil for which it is possible to use the $R$ matrix approach for generic $N$. This was shown in \cite{EGGKPSS} but we include it here for completeness. Following the labeled braid diagram in Figure \ref{fig: Braid Diagrams 31, 51, 819} the $R$ matrices we need are (Ignoring prefactors)
	\[
		R_{\textbf{0}, \textbf{b}}^{\textbf{b}, \textbf{0}} = x^{\frac{-b_1}{2}}, \qquad\qquad R_{\textbf{b}, \textbf{0}}^{\textbf{b}, \textbf{0}} = (-1)^{b_1} q^{\frac{b_1^2 - b_1}{2}} x^{-b_1}(x; q^{-1})_{b_1}, \qquad\qquad R_{\textbf{b}, \textbf{0}}^{\textbf{0}, \textbf{b}} = x^{\frac{-b_1}{2}}.
	\]
	Taking the quantum trace we get the general formula:
	\[
		F^{\G{sl}_N}_{3_1^r}(x, q) = x^{1 - N}\sum_{\textbf{b}} q^{-|\textbf{b}|} (-1)^{b_1} q^{\frac{1}{2}(b_1 - 1)b_1}  x^{-2b_1} (x; q^{-1})_{b_1}.
	\]
    Apply Weyl symmetry to send $x^{-1} \to a x$ and compute the sum over the $b_i$ variables for $i > 1$ to get:
	\[
		F^{\G{sl}_N}_{3_1^r}(x, q) = q^{(\log_q(x) + N)(N - 1)} \sum_{b_1 = 0}^{\infty} (-1)^{b_1} q^{\frac{1}{2}(b_1 + 1)b_1} q^{N b_1} x^{2b_1} \frac{(q^{b_1 + 1})_{N-2}(q^{-N} x^{-1}; q^{-1})_{b_1}}{(q)_{N-2}}.
	\]
	While this is already an expression which will work for all integer $N$, it's not quite in the right form for an $a$ deformation due to where the $N$'s appear in the $q$-Pochammers. We can fix this by observing that
	\[
	    \frac{(q^{b_1 + 1})_{N-2}}{(q)_{N-2}} = \frac{(q)_{N-2 + b_1}}{(q)_{b_1}(q)_{N-2}} = \frac{(q^{N - 1})_{b_1}}{(q)_{b_1}}
	\]
	which allows us to produce an $a$ deformed series
	\[
	    F_{3_1^r}(x, a, q) = q^{(\log_q(x) + \log_q(a))(\log_q(a) - 1)} \sum_{b_1 = 0}^{\infty} q^{b_1} x^{b_1} \frac{(a q^{-1})_{b_1}(a x, q)_{b_1}}{(q)_{b_1}}.
	\]
	This agrees with the $F_{3_1^l}(x, a, q)$ given\footnote{Note that \cite{EGGKPS} deals with the \textbf{left} handed trefoil whereas here we deal with the \textbf{right} handed one. Thus to relate this expression to \cite{EGGKPS}'s we need to send $(x, a, q) \to (x^{-1}, a^{-1}, q^{-1})$ and then apply Weyl symmetry. Additionally note that in \cite{EGGKPS} the prefactor has been dropped.} in \cite{EGGKPS} Section 5.2.
	
    
    For more complicated knots, we can make progress only after specialising $N$. For the $5_1$ knot, Figure \ref{fig: Braid Diagrams 31, 51, 819} shows
    \[
		F^{\G{sl}_N}_{5_1^r}(x, a = q^N, q) = x^{2(1 - N)} \sum_{\textbf{b}, \textbf{c}_1 ,\textbf{c}_2} q^{-|\textbf{b}|} x^{-b_1} R_{\textbf{b}, \textbf{0}}^{\textbf{c}_1, \textbf{b} - \textbf{c}_1} R_{\textbf{c}_1, \textbf{b} - \textbf{c}_1}^{\textbf{b} - \textbf{c}_2, \textbf{c}_2} R_{\textbf{b} - \textbf{c}_2, \textbf{c}_2}^{\textbf{0}, \textbf{b}}
	\]
	There are $2$ types of bounds on these summations. Internal bounds come from ensuring that labels are valid whereas external bounds come from crossings and the requirement that the label decreases along the overstrand. Recall that a label $\textbf{b}$ is valid if and only if $b_1 \geq b_2 \geq \cdots \geq b_{N - 1} \geq 0$ for all $i$. In this case, the labels are $\textbf{b}, \textbf{c}_1, \textbf{c}_2, \textbf{b} - \textbf{c}_1, \textbf{b} - \textbf{c}_2$ meaning that we need:
	\begin{align*}
	    0 & \leq b_{i + 1} \leq b_i \\
	    0, b_{i + 1} + c_{1, i} - b_i & \leq c_{1, i + 1} \leq c_{1, i}, b_{i + 1} \\
	    0, b_{i + 1} + c_{2, i} - b_i & \leq c_{2, i + 1} \leq c_{2, i}, b_{i + 1}  \\\
	\end{align*}
	In this case, the only external bound is $\textbf{c}_1 \geq \textbf{c}_2$. Using these bounds we compute $F^{\G{sl}_N}_{5_1^r}$ for $N = 2, 3, 4$ with the result shown in Table \ref{tab: 5 1} (Dropping the prefactor and applying Weyl symmetry). It is certainly possible to compute $F^N_K$ for larger $N$ but becomes increasingly time consuming as the number of summations is roughly quadratic in $N$. 
	
	    \begin{table}
        \centering
        \renewcommand{\arraystretch}{1.5}
        \begin{tabular}{ c | m{15cm} }
         & \multicolumn{1}{c}{$F_{5_1^r}(x, q)$} \\ \hline
        \shortstack{$\G{sl}_2$} & $1 + q x - (-1 + q) q^2 x^2 - (-1 + q) q^3 x^3 - (-1 + q) q^4 x^4 - 
                                    q^5 (-1 + q + q^4) x^5 - q^6 (-1 + q + q^4) x^6 + 
                                    q^7 (1 - q - q^4 + q^7) x^7 + q^8 (1 - q - q^4 + q^7) x^8$ \\  \hline
        \shortstack{$\G{sl}_3$} & $1 + q (1 + q) x + (q^2 + q^3 - q^5) x^2 - 
                                    q^3 (-1 - q + q^3 + q^4) x^3 - q^4 (-1 - q + q^3 + q^4) x^4 - 
                                    q^5 (-1 - q + q^3 + q^4 + q^8 + q^9) x^5 - 
                                    q^6 (-1 - q + q^3 + q^4 + q^8 + 2 q^9 + q^{10}) x^6 + 
                                    q^7 (1 + q - q^3 - q^4 - q^8 - 2 q^9 - 2 q^{10} + q^{12} + q^{13} + q^{14}) x^7 + 
                                    q^8 (1 + q - q^3 - q^4 - q^8 - 2 q^9 - 2 q^{10} - q^{11} + q^{12} + 2 q^{13} + 2 q^{14} + q^{15}) x^8$ \\ \hline
        \shortstack{$\G{sl}_4$} & $1 + q (1 + q + q^2) x + (q^2 + q^3 + 2 q^4 - q^7) x^2 + 
                                    q^3 (1 + q + 2 q^2 + q^3 - 2 q^5 - q^6 - q^7) x^3 + 
                                    q^4 (1 + q + 2 q^2 + q^3 + q^4 - 2 q^5 - 2 q^6 - 2 q^7) x^4 - 
                                    q^5 (-1 - q - 2 q^2 - q^3 - q^4 + q^5 + 2 q^6 + 3 q^7 + q^8 - q^{10} + q^{12} + q^{13} + q^{14}) x^5 - 
                                    q^6 (-1 - q - 2 q^2 - q^3 - q^4 + q^5 + q^6 + 3 q^7 + 2 q^8 + q^9 - q^{10} - q^{11} + q^{12} + 2 q^{13} + 3 q^{14} + 2 q^{15} + q^{16}) x^6 + 
                                    q^7 (1 + q + 2 q^2 + q^3 + q^4 - q^5 - q^6 - 2 q^7 - 2 q^8 - 2 q^9 + q^{11} - 2 q^{13} - 4 q^{14} - 4 q^{15} - 3 q^{16} + q^{18} + 2 q^{19} + q^{20} + q^{21}) x^7 + 
                                    q^8 (1 + q + 2 q^2 + q^3 + q^4 - q^5 - q^6 - 2 q^7 - q^8 - 2 q^9 - q^{10} - q^{13} - 4 q^{14} - 5 q^{15} - 5 q^{16} - 2 q^{17} + q^{18} + 4 q^{19} + 4 q^{20} + 4 q^{21} + 2 q^{22} + q^{23}) x^8$ \\  \hline 
                                    \shortstack{$\G{sl}_N$} & \[\sum_{d1, d2} a^{2d_2} x^{d_1 + 3d_2}q^{d_1 + d_2^2}\frac{(a^{-1} q, q^{-1})_{d_1 + d_2}(ax, q)_{d_1 + d_2}}{(q, q)_{d_1}(q, q)_{d_2}}\]
        \end{tabular}
        \caption{Computation of $F^N_{5_1^r}$ for small $N$. The $a$ deformation comes from \cite{EGGKPS} with orientation swapped.}
        \label{tab: 5 1}
    \end{table}
	
	We play a similar game for the $8_{19}$ knot. Note that there is a choice in the braid we use. In particular, the $8_{19}$ knot is usually represented as the closure of $\sigma_1^3\sigma_2\sigma_1^3\sigma_2$ where as we use the closure of $\sigma_1\sigma_2\sigma_1^2\sigma_2\sigma_1\sigma_2\sigma_1$. Whilst the first is conceptually simpler and cleaner, the second more naturally lends itself to these computations as more crossings are fixed by \textbf{0} labels and so we end up with fewer internal summations. Figure \ref{fig: Braid Diagrams 31, 51, 819} shows
    \[
		F^{\G{sl}_N}_{8_{19}}(x, a = q^N, q) = x^{3(N - 1)}\sum_{\textbf{a}, \textbf{b}, \textbf{c}_1 ,\textbf{c}_2} x^{-a_1 - b_1}q^{-|\textbf{a}| - |\textbf{b}|} R_{\textbf{a}, \textbf{b}, }^{\textbf{b} + \textbf{c}_1, \textbf{a} - \textbf{c}_1} R_{\textbf{b} + \textbf{c}_1, \textbf{a} - \textbf{c}_1}^{\textbf{b} + \textbf{c}_2, \textbf{a} - \textbf{c}_2} R_{\textbf{b} + \textbf{c}_2, \textbf{0}}^{\textbf{b} + \textbf{c}_2, \textbf{0}} R_{\textbf{b} + \textbf{c}_2, \textbf{a} - \textbf{c}_2}^{\textbf{a}, \textbf{b}}
	\]
    From here, similarly to the $5_1$ case, we compute $F^{\G{sl}_N}_{8_{19}}$ for $N = 2, 3, 4$ and results are shown in Table \ref{tab: 8 19}.
    
    One immediate observation we can make is that, even though we are only looking at a few values of $N$, we can already make prediction as to what an $a$ deformation should look like and we can check that these agree with the deformations presented in \cite{EGGKPS, EGGKPSS}.

    \begin{table}
        \centering
        \renewcommand{\arraystretch}{1.5}
        \begin{tabular}{ c | m{15cm} }
         & \multicolumn{1}{c}{$F_{8_{19}}(x, q)$} \\ \hline
        \shortstack{$\G{sl}_2$} & $1 + q x + q^2 x^2 + (q^3 - q^5) x^3 - q^4 (-1 + q^2 + q^3) x^4
                                    - q^5 (-1 + q^2 + q^3) x^5 - q^6 (-1 + q^2 + q^3) x^6 + q^7 (1 - q^2 - q^3 + q^7) x^7
                                    + q^8 (1 - q^2 - q^3 + q^7) x^8$ \\  \hline
        \shortstack{$\G{sl}_3$} & $1 + q (1 + q) x + q^2 (1 + q + q^2) x^2 - 
                                     q^3 (1 + q) (-1 - q^2 + q^4) x^3 - 
                                     q^4 (-1 - q - q^2 - q^3 + 2 q^5 + 2 q^6 + q^7) x^4 - 
                                     q^5 (1 + q) (-1 - q^2 + q^5 + 2 q^6 + q^7) x^5 - 
                                     q^6 (1 + q) (-1 - q^2 + q^5 + q^6 + 3 q^7) x^6 + 
                                     q^7 (1 + q)^2 (1 - q + 2 q^2 - 2 q^3 + 2 q^4 - 3 q^5 + 2 q^6 - 
                                        4 q^7 + 2 q^8 - 2 q^9 + 2 q^{10} - q^{11} + q^{12}) x^7 + 
                                     q^8 (1 + q + q^2 + q^3 - q^5 - 2 q^6 - 3 q^7 - 3 q^8 - 3 q^9 - 
                                        2 q^{10} + q^{11} + 2 q^{12} + 3 q^{13} + 2 q^{14} + q^{15}) x^8$ \\ \hline
        \shortstack{$\G{sl}_4$} & $1 + q (1 + q + q^2) x + q^2 (1 + q^2) (1 + q + q^2) x^2 - 
                                     q^3 (-1 - q - 2 q^2 - 2 q^3 - 2 q^4 - q^5 + q^7 + q^8) x^3 - 
                                     q^4 (1 + q + q^2) (-1 - q^2 - q^3 - q^4 + q^7 + q^8 + q^9) x^4 - 
                                     q^5 (1 + q + q^2) (-1 - q^2 - q^3 - q^4 - q^5 + q^7 + q^8 + 2 q^9 + 
                                        2 q^{10} + q^{11}) x^5 - 
                                     q^6 (-1 - q - 2 q^2 - 2 q^3 - 3 q^4 - 3 q^5 - 3 q^6 - q^7 + q^8 + 
                                        4 q^9 + 6 q^{10} + 8 q^{11} + 7 q^{12} + 5 q^{13} + 2 q^{14}) x^6 + 
                                     q^7 (1 + q + q^2) (1 + q^2 + q^3 + q^4 + q^5 + q^6 - q^8 - 2 q^9 - 
                                        3 q^{10} - 4 q^{11} - 3 q^{12} - 2 q^{13} - q^{14} + q^{16} + q^{17} + 
                                        q^{19}) x^7 + 
                                     q^8 (1 + q + q^2) (1 + q^2 + q^3 + q^4 + q^5 + q^6 - 2 q^9 - 
                                        3 q^{10} - 4 q^{11} - 4 q^{12} - 4 q^{13} - 2 q^{14} - q^{15} + 2 q^{17} + 
                                        2 q^{18} + 3 q^{19} + q^{20} + q^{21}) x^8$ \\  \hline 
        \shortstack{$\G{sl}_N$} & 
        $ \begin{aligned}
            1 & + \frac{(q^{-1} a)_1}{(q)_1} qx + \frac{(q^{-1} a)_2}{(q)_2}q^2x^2 + \left(\frac{(q^{-1} a)_3}{(q)_3}q^3 - a^2\frac{(q^{-1} a)_1}{(q)_1}q \right)x^3
            \\ & + \left(\frac{(q^{-1} a)_4}{(q)_4}q^4 - a^2 \frac{1-q^2}{1-q}\frac{(q^{-1} a)_2}{(q)_2}q^2 \right)x^4
            \\ & + \left(\frac{(q^{-1} a)_5}{(q)_5}q^5 - a^2 \frac{1-q^3}{1-q}\frac{(q^{-1} a)_3}{(q)_3}q^3 + a^4 \frac{(q^{-1} a)_1}{(q)_1}q\right)x^5
            \\ & + \left(\frac{(q^{-1} a)_6}{(q)_6}q^6 - a^2 \frac{1-q^4}{1-q}\frac{(q^{-1} a)_4}{(q)_4}q^4 + a^4(1 + 2q) \frac{(q^{-1} a)_2}{(q)_2}q^2 - a^5\frac{(q^{-1} a)_1}{(q)_1}q\right)x^6
        \end{aligned} $
        
        
        \end{tabular}
        \caption{Computation of $F^N_{8_{19}}$  for small $N$. The guess for the $a$ deformation is written suggestively to show the structure present.}
        \label{tab: 8 19}
    \end{table}
	
	\subsubsection{Non Torus Knots}
	
	Let us  now consider a couple of positive braid non torus knots. Up to $10$ crossing there are only $2$ of these namely, $10_{139}$ and $10_{152}$ and labeled braid diagrams for them are shown in Figure \ref{fig: Braid Diagrams 10139, 10152}. For each $N = 2, 3, 4$, we can compute $F^N_K$ as described above and the results are given in Tables \ref{tab: FK 10 139} and \ref{tab: FK 10 152}.
	
	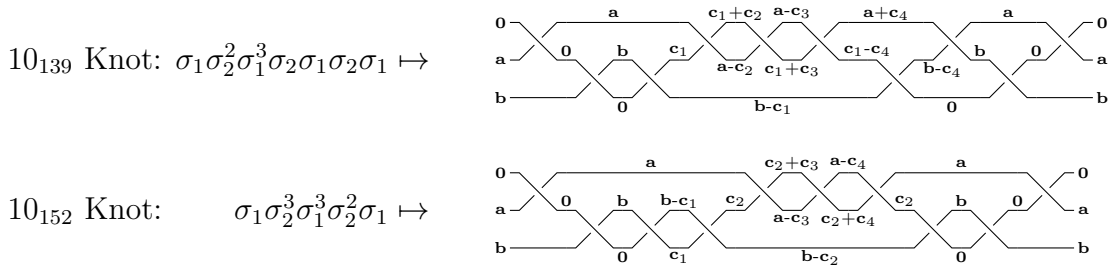
\begin{figure}[htp]
    	\centering
    	
    	\begin{tikzpicture}
    	
    	    \def\fsizeone{6}
    	    \def\fsizetwo{7}
    	    \def\size{0.5}
    	    \def\xshft{-2}
    	    
    	    
    	    \node at (-4 + \xshft, 0) {$10_{139}$ Knot: $\sigma_1\sigma_2^2\sigma_1^3\sigma_2\sigma_1\sigma_2\sigma_1 \mapsto$};
    	    
    	    \node at (\xshft - 0.25, \size) {\fontsize{\fsizeone}{\fsizetwo}\selectfont \textbf{0}};
    	    \node at (\xshft - 0.25, 0) {\fontsize{\fsizeone}{\fsizetwo}\selectfont \textbf{a}};
    	    \node at (\xshft - 0.25, -\size) {\fontsize{\fsizeone}{\fsizetwo}\selectfont \textbf{b}}; 
    	    
    	    \draw (\xshft - \size/4, -\size) -- (\xshft + 5*\size/4, -\size);
    	    
    	    \HoroUnderCross[\xshft, 0](\size)
    	    
    	    \node at (\xshft + 5*\size/4, \size/4) {\fontsize{\fsizeone}{\fsizetwo}\selectfont \textbf{0}};
    	    
    	    \draw (\xshft + 5*\size/4, \size) -- (\xshft + 17*\size/4, \size);
    	    \node at (\xshft + 10*\size/4, 5*\size/4) {\fontsize{\fsizeone}{\fsizetwo}\selectfont \textbf{a}};
    	    
    	    \HoroUnderCross[\xshft + 3*\size/2, -\size](\size)
    	    
    	    \node at (\xshft + 11*\size/4, -5*\size/4) {\fontsize{\fsizeone}{\fsizetwo}\selectfont \textbf{0}};
    	    
    	    \node at (\xshft + 11*\size/4, \size/4) {\fontsize{\fsizeone}{\fsizetwo}\selectfont \textbf{b}};
    	    
    	    \HoroUnderCross[\xshft + 6*\size/2, -\size](\size)
    	    
    	    \draw (\xshft + 17*\size/4, -\size) -- (\xshft + 37*\size/4, -\size);
    	    
    	    \node at (\xshft + 27*\size/4, -5*\size/4) {\fontsize{\fsizeone}{\fsizetwo}\selectfont \textbf{b}-\textbf{c}\textsubscript{1}};
    	    \node at (\xshft + 17*\size/4, \size/4) {\fontsize{\fsizeone}{\fsizetwo}\selectfont \textbf{c}\textsubscript{1}};
    	    
    	    \HoroUnderCross[\xshft + 9*\size/2, 0](\size)
    	    
    	    \node at (\xshft + 23*\size/4, 5*\size/4) {\fontsize{\fsizeone}{\fsizetwo}\selectfont \textbf{c}\textsubscript{1}+\textbf{c}\textsubscript{2}};
    	    \node at (\xshft + 23*\size/4, -\size/4) {\fontsize{\fsizeone}{\fsizetwo}\selectfont \textbf{a}-\textbf{c}\textsubscript{2}};
    	    
    	    \HoroUnderCross[\xshft + 12*\size/2, 0](\size)
    	    
    	    \node at (\xshft + 29*\size/4, 5*\size/4) {\fontsize{\fsizeone}{\fsizetwo}\selectfont \textbf{a}-\textbf{c}\textsubscript{3}};
    	    \node at (\xshft + 29*\size/4, -\size/4) {\fontsize{\fsizeone}{\fsizetwo}\selectfont \textbf{c}\textsubscript{1}+\textbf{c}\textsubscript{3}};
    	    
    	    \HoroUnderCross[\xshft + 15*\size/2, 0](\size)
    	    
    	    \draw (\xshft + 35*\size/4, \size) -- (\xshft + 43*\size/4, \size);
    	    \node at (\xshft + 39*\size/4, 5*\size/4) {\fontsize{\fsizeone}{\fsizetwo}\selectfont \textbf{a}+\textbf{c}\textsubscript{4}};
    	    
    	    \draw (\xshft + 35*\size/4, 0) -- (\xshft + 37*\size/4, 0);
    	    \node at (\xshft + 37*\size/4, \size/4) {\fontsize{\fsizeone}{\fsizetwo}\selectfont \textbf{c}\textsubscript{1}-\textbf{c}\textsubscript{4}};
    	    
    	    \HoroUnderCross[\xshft + 19*\size/2, -\size](\size)
    	    
    	    \draw (\xshft + 43*\size/4, -\size) -- (\xshft + 49*\size/4, -\size);
    	    \node at (\xshft + 46*\size/4, -5*\size/4) {\fontsize{\fsizeone}{\fsizetwo}\selectfont \textbf{0}};
    	    
    	    \node at (\xshft + 45*\size/4, -\size/4) {\fontsize{\fsizeone}{\fsizetwo}\selectfont \textbf{b}-\textbf{c}\textsubscript{4}};
    	    
    	    \HoroUnderCross[\xshft + 22*\size/2, 0](\size)
    	    
    	    \draw (\xshft + 49*\size/4, \size) -- (\xshft + 55*\size/4, \size);
    	    \node at (\xshft + 52*\size/4, 5*\size/4) {\fontsize{\fsizeone}{\fsizetwo}\selectfont \textbf{a}};
    	    
    	    \node at (\xshft + 49*\size/4, \size/4) {\fontsize{\fsizeone}{\fsizetwo}\selectfont \textbf{b}};
    	    
    	    \HoroUnderCross[\xshft + 25*\size/2, -\size](\size)
    	    
    	    \draw (\xshft + 55*\size/4, -\size) -- (\xshft + 61*\size/4, -\size);
    	   
    	    \node at (\xshft + 55*\size/4, \size/4) {\fontsize{\fsizeone}{\fsizetwo}\selectfont \textbf{0}};
    	    
    	    \HoroUnderCross[\xshft + 28*\size/2, 0](\size)
    	    
    	    \node at (\xshft + 62*\size/4, \size) {\fontsize{\fsizeone}{\fsizetwo}\selectfont \textbf{0}};
    	    \node at (\xshft + 62*\size/4, 0) {\fontsize{\fsizeone}{\fsizetwo}\selectfont \textbf{a}};
    	    \node at (\xshft + 62*\size/4, -\size) {\fontsize{\fsizeone}{\fsizetwo}\selectfont \textbf{b}};
    	    
    	    
    	    \def\yshft{-2}
    	    
    	    \node at (-4 + \xshft, \yshft) {$10_{152}$ Knot: \hspace{17pt} $\sigma_1\sigma_2^3\sigma_1^3\sigma_2^2\sigma_1 \mapsto$};
    	    
    	    \node at (\xshft - 0.25, \yshft + \size) {\fontsize{\fsizeone}{\fsizetwo}\selectfont \textbf{0}};
    	    \node at (\xshft - 0.25, \yshft) {\fontsize{\fsizeone}{\fsizetwo}\selectfont \textbf{a}};
    	    \node at (\xshft - 0.25, \yshft -\size) {\fontsize{\fsizeone}{\fsizetwo}\selectfont \textbf{b}};
    	    
    	    \draw (\xshft - \size/4, \yshft - \size) -- (\xshft + 5*\size/4, \yshft -\size);
    	    
    	    \HoroUnderCross[\xshft, \yshft](\size)
    	    
    	    \draw (\xshft + 5*\size/4, \yshft + \size) -- (\xshft + 23*\size/4, \yshft + \size);
    	    \node at (\xshft + 14*\size/4, \yshft + 5*\size/4) {\fontsize{\fsizeone}{\fsizetwo}\selectfont \textbf{a}};
    	    
    	    \node at (\xshft + 5*\size/4, \yshft + \size/4) {\fontsize{\fsizeone}{\fsizetwo}\selectfont \textbf{0}};
    	    
    	    \HoroUnderCross[\xshft + 3*\size/2, \yshft-\size](\size)
    	    
    	    \node at (\xshft + 11*\size/4, \yshft + \size/4) {\fontsize{\fsizeone}{\fsizetwo}\selectfont \textbf{b}};
    	    \node at (\xshft + 11*\size/4, \yshft + -5*\size/4) {\fontsize{\fsizeone}{\fsizetwo}\selectfont \textbf{0}};
    	    
    	    \HoroUnderCross[\xshft + 6*\size/2, \yshft-\size](\size)
    	    
    	    \node at (\xshft + 17*\size/4, \yshft + \size/4) {\fontsize{\fsizeone}{\fsizetwo}\selectfont \textbf{b}-\textbf{c}\textsubscript{1}};
    	    \node at (\xshft + 17*\size/4, \yshft + -5*\size/4) {\fontsize{\fsizeone}{\fsizetwo}\selectfont \textbf{c}\textsubscript{1}};
    	    
    	    \HoroUnderCross[\xshft + 9*\size/2, \yshft-\size](\size)
    	    
    	    \node at (\xshft + 23*\size/4, \yshft + \size/4) {\fontsize{\fsizeone}{\fsizetwo}\selectfont \textbf{c}\textsubscript{2}};
    	    
    	    \draw (\xshft + 23*\size/4, \yshft + -\size) -- (\xshft + 41*\size/4, \yshft + -\size);
    	    \node at (\xshft + 32*\size/4, \yshft + -5*\size/4) {\fontsize{\fsizeone}{\fsizetwo}\selectfont \textbf{b}-\textbf{c}\textsubscript{2}};
    	    
    	    \HoroUnderCross[\xshft + 12*\size/2, \yshft](\size)
    	    
    	    \node at (\xshft + 29*\size/4, \yshft + 5*\size/4) {\fontsize{\fsizeone}{\fsizetwo}\selectfont \textbf{c}\textsubscript{2}+\textbf{c}\textsubscript{3}};
    	    \node at (\xshft + 29*\size/4, \yshft + -\size/4) {\fontsize{\fsizeone}{\fsizetwo}\selectfont \textbf{a}-\textbf{c}\textsubscript{3}};
    	    
    	    \HoroUnderCross[\xshft + 15*\size/2, \yshft](\size)
    	    
    	    \node at (\xshft + 35*\size/4, \yshft + 5*\size/4) {\fontsize{\fsizeone}{\fsizetwo}\selectfont \textbf{a}-\textbf{c}\textsubscript{4}};
    	    \node at (\xshft + 35*\size/4, \yshft + -\size/4) {\fontsize{\fsizeone}{\fsizetwo}\selectfont \textbf{c}\textsubscript{2}+\textbf{c}\textsubscript{4}};
    	    
    	    \HoroUnderCross[\xshft + 18*\size/2, \yshft](\size)
    	    
    	    \draw (\xshft + 41*\size/4, \yshft + \size) -- (\xshft + 53*\size/4, \yshft + \size);
    	    \node at (\xshft + 47*\size/4, \yshft + 5*\size/4) {\fontsize{\fsizeone}{\fsizetwo}\selectfont \textbf{a}};
    	    
    	    \node at (\xshft + 41*\size/4, \yshft + \size/4) {\fontsize{\fsizeone}{\fsizetwo}\selectfont \textbf{c}\textsubscript{2}};
    	    
    	    \HoroUnderCross[\xshft + 21*\size/2, \yshft-\size](\size)
    	    
    	    \node at (\xshft + 47*\size/4, \yshft + \size/4) {\fontsize{\fsizeone}{\fsizetwo}\selectfont \textbf{b}};
    	    \node at (\xshft + 47*\size/4, \yshft + -5*\size/4) {\fontsize{\fsizeone}{\fsizetwo}\selectfont \textbf{0}};
    	    
    	    \HoroUnderCross[\xshft + 24*\size/2, \yshft-\size](\size)
    	    
    	    \node at (\xshft + 53*\size/4, \yshft + \size/4) {\fontsize{\fsizeone}{\fsizetwo}\selectfont \textbf{0}};
    	    
    	    \draw (\xshft + 53*\size/4, \yshft - \size) -- (\xshft + 59*\size/4, \yshft - \size);
    	    
    	    \HoroUnderCross[\xshft + 27*\size/2, \yshft](\size)
    	    
    	    \node at (\xshft + 60*\size/4, \yshft + \size) {\fontsize{\fsizeone}{\fsizetwo}\selectfont \textbf{0}};
    	    \node at (\xshft + 60*\size/4, \yshft) {\fontsize{\fsizeone}{\fsizetwo}\selectfont \textbf{a}};
    	    \node at (\xshft + 60*\size/4, \yshft -\size) {\fontsize{\fsizeone}{\fsizetwo}\selectfont \textbf{b}};

        \end{tikzpicture}
        \captionsetup{singlelinecheck=off}
    	\caption[.]{Labeled Braid diagrams for the $10_{139}$ and $10_{152}$ knots.}
    	\label{fig: Braid Diagrams 10139, 10152}
    \end{figure}
	
	We are interested in trying to predict the $a$ deformation from the computations at small $N$. Note that as stated so far this problem is inherently underdetermined as we only know the solution for $N = 1, 2, 3, 4$ and so any guess can always be modified by a term containing $(q - a)(q^2 - a)(q^3 - a)(q^4 - a)$ (Or similar expressions) which are exactly $0$ in these cases. However we do have one extra constraint, namely the $a = q^N, q \to 1$ limit for each $N$ upon which the series should collapse to $\frac{1}{\Delta_K(X)^{N - 1}}$. Additionally, making use of the Knot-Quivers correspondence \cite{KRSS1, KRSS2, Kuch, EGGKPSS} mentioned earlier, we can provide a useful ansatz \eqref{eq: Quiver Form} for this $a$ deformation. In practise, instead of using Equation \eqref{eq: Quiver Form}, a more useful ansatz seems to be (Ignoring the prefactor for a moment)
	\[
	    \sum_{\textbf{d}} (-q)^{\frac{1}{2}\textbf{d}M\textbf{d}^T}q^{\textbf{q}\cdot \textbf{d}}a^{\textbf{a}\cdot \textbf{d}}x^{\textbf{x}\cdot \textbf{d}} \frac{(q^{-1}a;q)_{\textbf{d}}}{(q; q)_{\textbf{d}}}.
	\]
	This is essentially equivalent to Equation \eqref{eq: Quiver Form} but hard-codes the $a = q$ specialisation as $F_K(x, q, q) = 1$. As we can see in Tables \ref{tab: 8 19}, \ref{tab: FK 10 139}, \ref{tab: FK 10 152} this simple anstaz does an excellent job of allowing us to guess the $a$ deformations of the first couple of terms. We can check that the $a$ deformation of $8_{19}$ is correct and we conjecture that the $a$ deformations for the $10_{139}, 10_{152}$ will also yield the correct series. Unfortunately it is difficult to pass from these initial terms to a full Quiver form so we only have these $a$ deformations perturbatively.
	
	\begin{table}
        \centering
        \renewcommand{\arraystretch}{1.5}
        \begin{tabular}{ c | m{15cm} }
         & \multicolumn{1}{c}{$F_{10_{139}}(x, q)$} \\ \hline
        \shortstack{$\G{sl}_2$} & $1 + q x + q^2 x^2 - q^3 (-1 + 2 q^2) x^3 + (-1 + q) q^4 (-1 - q + q^2) x^4 - q^5 (-1 + 2 q^2 - q^3 + q^4) x^5 + q^6 (1 - 2 q^2 + q^3 - q^4 + 2 q^5 + q^6) x^6 - q^7 (-1 + 2 q^2 - q^3 + q^4 - 2 q^5 + q^6 + 2 q^7) x^7 + q^8 (1 - 2 q^2 + q^3 - q^4 + 2 q^5 - q^6 - 2 q^7 + 3 q^8 + q^9) x^8$ \\  \hline
        \shortstack{$\G{sl}_3$} & $1 + q (1 + q) x + q^2 (1 + q + q^2) x^2 - (-1 + q) q^3 (1 + q)^2 (1 + 2 q^2) x^3 + q^4 (1 + q + q^2 + q^3 - q^4 - 4 q^5 - q^6 + q^7) x^4 - q^5 (1 + q)^2 (-1 + q - 2 q^2 + 2 q^3 - q^4 + 3 q^5 - 2 q^6 + q^7) x^5 + q^6 (1 + q + q^2 + q^3 - q^4 - 3 q^5 - 2 q^6 - 2 q^7 - q^8 + 2 q^9 + 2 q^{10} + 3 q^{11} + q^{12}) x^6 - q^7 (1 + q) (-1 - q^2 + q^4 + 2 q^5 + q^7 + 2 q^8 - 3 q^9 - 2 q^{10} - q^{11} + q^{12} + 2 q^{13}) x^7 + q^8 (1 + q) (1 + q^2 - q^4 - 2 q^5 - q^7 - q^8 + 4 q^{10} + 2 q^{11} - 2 q^{12} - q^{13} - q^{14} + 3 q^{15} + q^{16}) x^8$ \\ \hline
        \shortstack{$\G{sl}_4$} & $1 + q (1 + q + q^2) x + q^2 (1 + q^2) (1 + q + q^2) x^2 - q^3 (1 + q) (-1 - 2 q^2 - 2 q^4 + q^5 + 2 q^7) x^3 + (-1 + q) q^4 (1 + q) (1 + q + q^2) (-1 - 2 q^2 - q^3 - 3 q^4 - q^5 - 2 q^6 + q^7) x^4 - q^5 (1 + q + q^2) (-1 - q^2 - q^3 - q^4 - q^5 + q^6 + 2 q^7 + 3 q^8 + q^9 + q^{10} - q^{11} + q^{12}) x^5 + q^6 (1 + q^2) (1 + q + q^2 + q^3 + 2 q^4 + 2 q^5 - 3 q^7 - 5 q^8 - 4 q^9 - 3 q^{10} - q^{11} - q^{12} + 2 q^{13} + 2 q^{14} + 3 q^{15} + q^{16}) x^6 - q^7 (1 + q + q^2) (-1 - q^2 - q^3 - q^4 - q^5 + q^7 + 3 q^8 + 2 q^9 + 4 q^{10} + 2 q^{11} + 2 q^{12} - 2 q^{13} - q^{14} - 6 q^{15} - q^{16} - q^{17} + q^{18} + 2 q^{19}) x^7 + q^8 (1 + q + q^2) (1 + q^2 + q^3 + q^4 + q^5 - q^7 - 2 q^8 - 2 q^9 - 4 q^{10} - 3 q^{11} - 4 q^{12} + q^{13} + q^{14} + 7 q^{15} + 3 q^{16} + 5 q^{17} - 2 q^{18} - 2 q^{19} - q^{21} + 3 q^{22} + q^{23}) x^8$ \\ \hline
        \shortstack{$\G{sl}_N$} &
        $ \begin{aligned}
            1 & + \frac{(q^{-1} a)_1}{(q)_1} qx + \frac{(q^{-1} a)_2}{(q)_2}q^2x^2 + \left(\frac{(q^{-1} a)_3}{(q)_3} q^3 - 2a^2\frac{(q^{-1} a)_1}{(q)_1} q\right)x^3
            \\ & + \left(\frac{(q^{-1} a)_4}{(q)_4} q^4 - 2a^2\frac{1 - q^2}{1 - q}\frac{(q^{-1} a)_2}{(q)_2} q^2 + 3a^3\frac{(q^{-1} a)_1}{(q)_1} q\right)x^4 {\quad \quad \quad}
            \\ & + \left(\frac{(q^{-1} a)_5}{(q)_5} q^5 - 2a^2\frac{1 - q^3}{1 - q}\frac{(q^{-1} a)_3}{(q)_3} q^3 + 3a^3\frac{1 - q^2}{1 - q}\frac{(q^{-1} a)_2}{(q)_2} q^2- 2a^4\frac{(q^{-1} a)_1}{(q)_1} q\right)x^5
            \\ & + \Bigg(\frac{(q^{-1} a)_6}{(q)_6} q^6 - 2a^2\frac{1 - q^4}{1 - q}\frac{(q^{-1} a)_4}{(q)_4} q^4 + 3a^3\frac{1 - q^3}{1 - q}\frac{(q^{-1} a)_3}{(q)_3} q^2
            \\ & \hspace{2cm} - a^4(2 - q - q^2)\frac{(q^{-1} a)_1}{(q)_1} q\Bigg)x^6
        \end{aligned} $
        
        \end{tabular}
        \caption{$F_K^N$ invariant for the $10_{139}$ knot with lie algebra $\G{sl}_N$ for small $N$.}
        \label{tab: FK 10 139}
    \end{table}

    \begin{table}
        \centering
        \renewcommand{\arraystretch}{1.5}
        \begin{tabular}{ c | m{15cm} }
         & \multicolumn{1}{c}{$F_{10_{152}}(x, q)$} \\ \hline
        \shortstack{$\G{sl}_2$} & $1 + q x + q^2 (1 + q) x^2 - q^3 (-1 - q + 3 q^2) x^3 + q^4 (1 + q - 2 q^2 + 2 q^3) x^4
                                    - q^5 (-1 - q + 2 q^2 + 5 q^4 + q^5) x^5 + q^6 (1 + q - 2 q^2 + q^3 - 4 q^4 + 6 q^5 + 4 q^6) x^6
                                    - q^7 (-1 - q + 2 q^2 - q^3 + 6 q^4 - 4 q^5 + 5 q^6 + 9 q^7 + 3 q^8) x^7
                                    + q^8 (1 + q - 2 q^2 + q^3 - 5 q^4 + 5 q^5 - q^6 + 15 q^8 + 11 q^9 + 3 q^{10}) x^8$ \\  \hline
        \shortstack{$\G{sl}_3$} & $1 + q (1 + q) x + q^2 (1 + q + 2 q^2 + q^3) x^2 - q^3 (1 + q) (-1 - 2 q^2 - q^3 + 3 q^4) x^3 
                                    + q^4 (1 + q + 2 q^2 + 3 q^3 + q^4 - 4 q^5 + 2 q^7) x^4 - q^5 (1 + q) (-1 - 2 q^2 - q^3 + q^5 + 2 q^6 + q^7 + 5 q^8 + q^9) x^5
                                    + q^6 (1 + q + 2 q^2 + 3 q^3 + q^4 - q^5 + q^6 - 5 q^7 - 9 q^8 - 2 q^9 + 5 q^{10} + 10 q^{11} + 4 q^{12}) x^6
                                    - q^7 (1 + q) (-1 - 2 q^2 - q^3 + q^5 - 2 q^6 + 3 q^7 + 10 q^8 - q^9 - q^{10} - q^{11} + 8 q^{12} + 9 q^{13} + 3 q^{14}) x^7
                                    + q^8 (1 + q + 2 q^2 + 3 q^3 + q^4 - q^5 + q^6 - q^7 - 8 q^8 - 12 q^9 - 3 q^{10} + 5 q^{11} - q^{12} - 3 q^{13} + 9 q^{14} + 26 q^{15} + 29 q^{16} + 14 q^{17} + 3 q^{18}) x^8$ \\ \hline
        \shortstack{$\G{sl}_4$} & $1 + q(1 + q + q^2) x + q^2(1 + q + q^2) (1 + q^2 + q^3) x^2 + q^3(1 + q + 2 q^2 + 3 q^3 + 4 q^4 + 4 q^5 - 2 q^7 - 3 q^8) x^3
                                    + q^4(1 + q + q^2) (1 + q^2 + 2 q^3 + 2 q^4 + 2 q^5 - 2 q^7 - 2 q^8 + 2 q^9) x^4
                                    - q^5(1 + q + q^2) (-1 - q^2 - 2 q^3 - 2 q^4 - 3 q^5 - q^6 + 2 q^8 + q^9 + 3 q^{10} + q^{11} + 5 q^{12} + q^{13}) x^5
                                    + q^6(1 + q + 2 q^2 + 3 q^3 + 5 q^4 + 7 q^5 + 7 q^6 + 6 q^7 + 3 q^8 + q^9 - 5 q^{10} - 9 q^{11} - 17 q^{12} - 12 q^{13} - 9 q^{14} + 7 q^{15} + 9 q^{16} + 10 q^{17} + 4 q^{18}) x^6 
                                    - q^7(1 + q + q^2) (-1 - q^2 - 2 q^3 - 2 q^4 - 3 q^5 - 2 q^6 - 2 q^7 - q^8 - 2 q^9 + 4 q^{10} + 5 q^{11} + 12 q^{12} + 5 q^{13} + 7 q^{14} - 8 q^{15} + 2 q^{16} + 2 q^{17} + 8 q^{18} + 9 q^{19} + 3 q^{20}) x^7
                                    + q^8(1 + q + q^2) (1 + q^2 + 2 q^3 + 2 q^4 + 3 q^5 + 2 q^6 + 2 q^7 + 2 q^8 + 3 q^9 - 2 q^{10} - 4 q^{11} - 12 q^{12} - 8 q^{13} - 12 q^{14} + 2 q^{15} - 4 q^{16} + 4 q^{17} - 6 q^{18} - 2 q^{19} + 9 q^{20} + 11 q^{21} + 18 q^{22} + 11 q^{23} + 3 q^{24}) x^8$ \\ \hline
            
            \shortstack{$\G{sl}_N$} & 
            
            $ \begin{aligned}
            1 & + \frac{(q^{-1} a)_1}{(q)_1} qx + \left(\frac{(q^{-1} a)_2}{(q)_2} q^2 + a\frac{(q^{-1} a)_1}{(q)_1} q \right)x^2
            \\ & + \left(\frac{(q^{-1} a)_3}{(q)_3} q^3 + a\frac{1 - q^2}{1 - q}\frac{(q^{-1} a)_2}{(q)_2} q^2 - 4a^2\frac{(q^{-1} a)_1}{(q)_1} q\right)x^3
            \\ & + \left(\frac{(q^{-1} a)_4}{(q)_4} q^4 + a\frac{1 - q^3}{1 - q}\frac{(q^{-1} a)_3}{(q)_3} q^3 - a^2(3 + 4q)\frac{(q^{-1} a)_2}{(q)_2} q^2 + 5a^3\frac{(q^{-1} a)_1}{(q)_1} q\right)x^4
            \\ & + \Bigg(\frac{(q^{-1} a)_5}{(q)_5} q^5 + a\frac{1 - q^4}{1 - q}\frac{(q^{-1} a)_4}{(q)_4} q^4 - 3a^2\frac{1 - q^3}{1 - q}\frac{(q^{-1} a)_3}{(q)_3} q^3
            \\ & \hspace{2cm} + a^3(2 + q - q^2)\frac{(q^{-1} a)_2}{(q)_2} q^2 - 4a^4\frac{(q^{-1} a)_1}{(q)_1} q\Bigg)x^5
        \end{aligned} $
        
        \end{tabular}
        \caption{$F_K$ invariant for the $10_{152}$ knot and symmetric series of representations on $\G{sl}_N$ for small $N$.}
        \label{tab: FK 10 152}
    \end{table}
    
    \subsection{Stratified state sum}
        
        An immediate question we encounter is what happens if we try to apply this technology to non positive braid knots. While in general this state sum will not converge, in certain cases can get conditional convergence by stratifying the summation \cite{Park2}.
        
        The key point is that the braid group representation on $V^{\otimes (n + 1)}$ is not irreducible. In particular if we define the total weight of a set of state by the sum of the labels $\textbf{w} = \textbf{b}_1 + \cdots + \textbf{b}_n$, we immediately observe that \textbf{w} is fixed by the action of the $R$, $R^{-1}$ matrices. Hence, letting $V^{\otimes (n + 1)}_{\textbf{w}}$ denote the subspace of total weight \textbf{w} we can stratify by total weight, yielding
        \[
            \wt{\Tr}^{q, x}_{V^{\otimes (n + 1)}}(\beta) = \lim_{\eta \to 1} \sum_{\textbf{w}} \eta^{|\textbf{w}|} \wt{\Tr}^q_{V^{\otimes n}_{\textbf{w}}}(\beta).
        \]
        When this converges we conjecture that it produces the correct $F_K$ invariant. Two examples of this are the $m(5_2)$ and $m(7_3)$ knots. For both of these knots, the Alexander polynomial is not monic which means that we should expect each coefficient in $x$ to itself be a power series in $q$. Using the standard notation
        \[
            F^N_K = x^{\frac{(N - 1)}{2}(\# Tr + \# R^{-1} - \# R)} \sum_{x = 0}^{\infty} f^N_k(n; q) x^n
        \]
        we find:
        \begin{align*}
            f^3_{m(5_2)}(0; q) & = 1 - q - q^2 + q^3 + q^4 + q^5 - q^6 - q^7 - q^8 - q^9 + q^{10} + q^{11} + q^{12} + q^{13} + O(q^{14})
            \\ f^3_{m(5_2)}(1; q) & = 2 q - 4 q^3 - q^4 + 3 q^5 + 5 q^6 + 2 q^7 - 2 q^8 - 4 q^9 - 6 q^{10} - 3 q^{11} + q^{12} + 3 q^{13} + O(q^{14})
            \\ f^3_{m(5_2)}(2; q) & = 3 q^2 + q^3 - 5 q^4 - 7 q^5 + q^6 + 11 q^7 + 11 q^8 + 5 q^9 - 6 q^{10} - 14 q^{11} - 16 q^{12} + O(q^{13})
        \end{align*}
        and 
        \begin{align*}
            f^3_{m(7_3)}(0; q) & = 1 - q - q^2 + q^3 + q^4 + q^5 - q^6 - q^7 - q^8 - q^9 + q^{10} + q^{11} + q^{12} + q^{13} + O(q^{14})
            \\ f^3_{m(7_3)}(1; q) & = 2q - q^2 - 3q^3 + 3q^5 + 4q^6 - 2q^8 - 4q^9 - 4q^{10} -q^{11} + 2q^{12} + 3q^{13} + 4q^{14} + O(q^{14})
            \\ f^3_{m(7_3)}(2; q) & = -q + 4q^2 - 5q^4 - 4q^5 + 2q^6 + 9q^7 + 5q^8 + 2q^9 - 6q^{10} - 10q^{11} + O(q^{12})
        \end{align*}
        Further terms can be most easily computed via recursion on the quantum A-polynomial as described in \cite{EGGKPS} with the terms here giving the initial conditions. In the $m(5_2)$ case, we can cross check this with the quiver form for $F_{m(5_2)}(x, a, q)$ given in \cite{EGGKPSS}. With a little effort\footnote{We take the quiver form, consider all terms contributing to a particular power of $x$ and apply the following identity:
        \[
            \sum_{d1, d2 = 0}^{\infty} \frac{\big(-q^{\frac{1}{2}}\big)^{d1^2 + 2 d1d2 + 2d2^2} b^{d1}a^{d2}q^{\frac{d1}{2}}}{(q, q)_{d1}(q, q)_{d2}} = \sum_{d = 0}^{\infty} (-1)^d q^{\frac{d(d + 1)}{2}}b^d \frac{(a b^{-1}, q)_d}{(q, q)_d}.
        \]}
        we simplify the form yielding the infinite sum expressions:
        \begin{align*}
            f_{m(5_2)}(0; a, q) & = \sum_{d = 0}^{\infty} (-1)^d q^{\frac{d(d + 1)}{2}}\frac{(a q^{-1}, q)_d}{(q, q)_d}
            \\ f_{m(5_2)}(1; a, q) & = \sum_{d = 0}^{\infty} (-1)^d q^{\frac{d(d + 1)}{2}}\frac{\big(q(2 - q^{-d})(a q^{-1}, q)_d - a(1 + a q^d)(a, q)_d\big)}{(q, q)_1(q, q)_d}
        \end{align*}
        It is an easy check that when $a = q$, $f_{5_2}(0; a = q, q) = 1$ and $f_{5_2}(1; a = q, q) = 0$. Similarly, specialising $a = q^2$ yields the expressions in \cite{Park2} and $a = q^3$ yields the expressions above.
        
        \begin{figure}[htp]
    	\centering
    	
    	\begin{tikzpicture}
    	
    	    \def\fsizeone{9}
    	    \def\fsizetwo{10}
    	    \def\size{0.75}
    	    \def\xshft{-2}
    	    
    	    
    	    \node at (-5 + \xshft, 0) {$m(5_2)$ Knot: $\sigma_1^2\sigma_2\sigma_1^{-1}\sigma_2\sigma_1 \mapsto$};
    	    
    	    \node at (\xshft - \size/2, \size) {\fontsize{\fsizeone}{\fsizetwo}\selectfont \textbf{0}};
    	    \node at (\xshft - \size/2, 0) {\fontsize{\fsizeone}{\fsizetwo}\selectfont \textbf{a}};
    	    \node at (\xshft - \size/2, -\size) {\fontsize{\fsizeone}{\fsizetwo}\selectfont \textbf{b}}; 
    	    
    	    \draw (\xshft - \size/4, -\size) -- (\xshft + 11*\size/4, -\size);
    	    
    	    \HoroUnderCross[\xshft, 0](\size)
    	    
    	    \node at (\xshft + 5*\size/4, 5*\size/4) {\fontsize{\fsizeone}{\fsizetwo}\selectfont \textbf{a}};
    	    \node at (\xshft + 5*\size/4, -\size/4) {\fontsize{\fsizeone}{\fsizetwo}\selectfont \textbf{0}};
    	    
    	    \HoroUnderCross[\xshft + 3*\size/2, 0](\size)
    	    
    	    \node at (\xshft + 11*\size/4, -\size/4) {\fontsize{\fsizeone}{\fsizetwo}\selectfont \textbf{c}\textsubscript{1}};
    	    
    	    \draw (\xshft + 11*\size/4, \size) -- (\xshft + 17*\size/4, \size);
    	    \node at (\xshft + 14*\size/4, 5*\size/4) {\fontsize{\fsizeone}{\fsizetwo}\selectfont \textbf{a}-\textbf{c}\textsubscript{1}};
    	    
    	    \HoroUnderCross[\xshft + 6*\size/2, -\size](\size)
    	    
    	    \node at (\xshft + 31*\size/8, \size/4) {\fontsize{\fsizeone}{\fsizetwo}\selectfont \textbf{b}+\textbf{c}\textsubscript{1}};
    	    
    	    \draw (\xshft + 17*\size/4, -\size) -- (\xshft + 23*\size/4, -\size);
    	    \node at (\xshft + 20*\size/4, -5*\size/4) {\fontsize{\fsizeone}{\fsizetwo}\selectfont \textbf{0}};
    	    
    	    \HoroOverCross[\xshft + 9*\size/2, 0](\size)
    	    
    	    \node at (\xshft + 23*\size/4, \size/4) {\fontsize{\fsizeone}{\fsizetwo}\selectfont \textbf{b}};
    	    
    	    \draw (\xshft + 23*\size/4, \size) -- (\xshft + 29*\size/4, \size);
    	    \node at (\xshft + 26*\size/4, 5*\size/4) {\fontsize{\fsizeone}{\fsizetwo}\selectfont \textbf{a}};

    	    \HoroUnderCross[\xshft + 12*\size/2, -\size](\size)
    	    
    	    \node at (\xshft + 29*\size/4, \size/4) {\fontsize{\fsizeone}{\fsizetwo}\selectfont \textbf{0}};
    	    
    	    \draw (\xshft + 29*\size/4, -\size) -- (\xshft + 35*\size/4, -\size);
    	    
    	    \HoroUnderCross[\xshft + 15*\size/2, 0](\size)
    	    
    	    \node at (\xshft + 36*\size/4, \size) {\fontsize{\fsizeone}{\fsizetwo}\selectfont \textbf{0}};
    	    \node at (\xshft + 36*\size/4, 0) {\fontsize{\fsizeone}{\fsizetwo}\selectfont \textbf{a}};
    	    \node at (\xshft + 36*\size/4, -\size) {\fontsize{\fsizeone}{\fsizetwo}\selectfont \textbf{b}};
    	    
    	    
    	    \def\yshft{-2.5}
    	    
    	    \node at (-5 + \xshft, \yshft) {$m(7_3)$ Knot: $\sigma_1^4\sigma_2\sigma_1^{-1}\sigma_2\sigma_1 \mapsto$};
    	    
    	    \node at (\xshft - 4*\size/2, \size + \yshft) {\fontsize{\fsizeone}{\fsizetwo}\selectfont \textbf{0}};
    	    \node at (\xshft - 4*\size/2, 0 + \yshft) {\fontsize{\fsizeone}{\fsizetwo}\selectfont \textbf{a}};
    	    \node at (\xshft - 4*\size/2, -\size + \yshft) {\fontsize{\fsizeone}{\fsizetwo}\selectfont \textbf{b}};
    	    
    	    \draw (\xshft - 7*\size/4, -\size + \yshft) -- (\xshft + 17*\size/4, -\size + \yshft);
    	    
    	    \HoroUnderCross[\xshft - 3*\size/2, 0 + \yshft](\size)
    	    
    	    \node at (\xshft - \size/4, 5*\size/4 + \yshft) {\fontsize{\fsizeone}{\fsizetwo}\selectfont \textbf{a}};
    	    \node at (\xshft - \size/4, -\size/4 + \yshft) {\fontsize{\fsizeone}{\fsizetwo}\selectfont \textbf{0}};
    	    
    	    \HoroUnderCross[\xshft, 0 + \yshft](\size)
    	    
    	    \node at (\xshft + 5*\size/4, 5*\size/4 + \yshft) {\fontsize{\fsizeone}{\fsizetwo}\selectfont \textbf{c}\textsubscript{1}};
    	    \node at (\xshft + 5*\size/4, -\size/4 + \yshft) {\fontsize{\fsizeone}{\fsizetwo}\selectfont \textbf{a}-\textbf{c}\textsubscript{1}};
    	    
    	    \HoroUnderCross[\xshft + 3*\size/2, 0 + \yshft](\size)
    	    
    	    \node at (\xshft + 11*\size/4, 5*\size/4 + \yshft) {\fontsize{\fsizeone}{\fsizetwo}\selectfont \textbf{a}-\textbf{c}\textsubscript{2}};
    	    \node at (\xshft + 11*\size/4, -\size/4 + \yshft) {\fontsize{\fsizeone}{\fsizetwo}\selectfont \textbf{c}\textsubscript{2}};
    	    
    	    \HoroUnderCross[\xshft + 6*\size/2, 0 + \yshft](\size)
    	    
    	    \draw (\xshft + 17*\size/4, \size + \yshft) -- (\xshft + 23*\size/4, \size + \yshft);
    	    \node at (\xshft + 20*\size/4, 5*\size/4 + \yshft) {\fontsize{\fsizeone}{\fsizetwo}\selectfont \textbf{a}-\textbf{c}\textsubscript{3}};
    	    
    	    \node at (\xshft + 17*\size/4, -\size/4 + \yshft) {\fontsize{\fsizeone}{\fsizetwo}\selectfont \textbf{c}\textsubscript{3}};

    	    \HoroUnderCross[\xshft + 9*\size/2, -\size + \yshft](\size)
    	    
    	    \node at (\xshft + 43*\size/8, \size/4 + \yshft) {\fontsize{\fsizeone}{\fsizetwo}\selectfont \textbf{b}+\textbf{c}\textsubscript{3}};
    	    
    	    \draw (\xshft + 23*\size/4, -\size + \yshft) -- (\xshft + 29*\size/4, -\size + \yshft);
    	    \node at (\xshft + 26*\size/4, -5*\size/4 + \yshft) {\fontsize{\fsizeone}{\fsizetwo}\selectfont \textbf{0}};
    	    
     	    \HoroOverCross[\xshft + 12*\size/2, 0 + \yshft](\size)
    	    
    	    \node at (\xshft + 29*\size/4, \size/4 + \yshft) {\fontsize{\fsizeone}{\fsizetwo}\selectfont \textbf{b}};
    	    
    	    \draw (\xshft + 29*\size/4, \size + \yshft) -- (\xshft + 35*\size/4, \size + \yshft);
    	    \node at (\xshft + 32*\size/4, 5*\size/4 + \yshft) {\fontsize{\fsizeone}{\fsizetwo}\selectfont \textbf{a}};

    	    \HoroUnderCross[\xshft + 15*\size/2, -\size + \yshft](\size)
    	    
    	    \node at (\xshft + 35*\size/4, \size/4 + \yshft) {\fontsize{\fsizeone}{\fsizetwo}\selectfont \textbf{0}};
    	    
    	    \draw (\xshft + 35*\size/4, -\size + \yshft) -- (\xshft + 41*\size/4, -\size + \yshft);
    	    
    	    \HoroUnderCross[\xshft + 18*\size/2, 0 + \yshft](\size)
    	    
    	    \node at (\xshft + 42*\size/4, \size + \yshft) {\fontsize{\fsizeone}{\fsizetwo}\selectfont \textbf{0}};
    	    \node at (\xshft + 42*\size/4, 0 + \yshft) {\fontsize{\fsizeone}{\fsizetwo}\selectfont \textbf{a}};
    	    \node at (\xshft + 42*\size/4, -\size + \yshft) {\fontsize{\fsizeone}{\fsizetwo}\selectfont \textbf{b}};

        \end{tikzpicture}
        \captionsetup{singlelinecheck=off}
    	\caption[.]{Labeled Braid diagrams for the $m(5_2)$ and $m(7_3)$ knots.}
    	\label{fig: Braid Diagrams m52, m73}
    \end{figure}
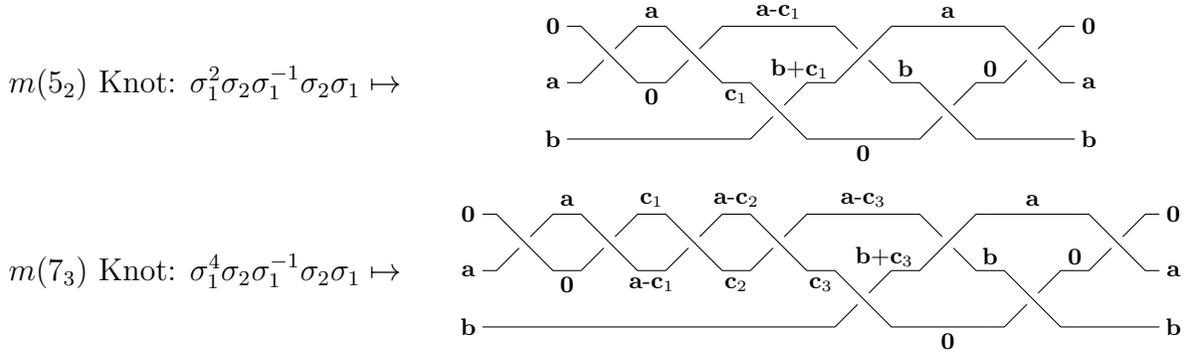

    In general, the stratified state sum will converge provided that, for each incoming state $\textbf{b}_i$, either\footnote{If both the minimal $x^{-1}$ and $q$ powers strictly increase as $\textbf{b}_i$ increase then the sum converges absolutely.} the minimal $x^{-1}$ or $q$ power strictly increases as $\textbf{b}_i$ increases. It would be nice to have an easy way to check this in a labelled braid diagram but unfortunately we can only easily check the behaviour of the $x$ power. 
    
\section{The Inverse State Sum Technique} \label{sec: inverse state sum}
    In cases where the stratified state sum does not converge, there is one final technique we can try to apply. This technique, called the Inverse State Sum was recently introduced in \cite{Park3} for the $\G{sl}_2$ case and extended the $R$ matrix computations to signed braid links. We take a very computational approach here as the proof in \cite{Park3} does not currently extend to $\G{sl}_N$ but we conjecture that the method still works.
    
    The basic observation is that a large number of $q$-Hypergeometric identities remain true if you invert the direction of the sum. The simplest example of this phenomenon is
    \[
        \sum_{n = 0}^{\infty} q^n = \frac{1}{1 - q} = -\frac{1}{q}\frac{1}{1 - \frac{1}{q}} = -\sum_{n < 0} q^{n}.
    \]
    This example can be easily generalized to:
    \begin{align*}
        \sum_{n = 0}^{\infty} \frac{(q^{n + m}, q^{-1})_m}{(q;q)_m} x^n & = \frac{1}{(x;q)_{m + 1}}
        \\ & = \frac{(-1)^{m + 1}}{x^{m + 1}q^{\frac{m(m + 1)}{2}}} \frac{1}{(x^{-1};q^{-1})_{m+1}}
        \\ & = \frac{(-1)^{m + 1}}{x^{m + 1}q^{\frac{m(m + 1)}{2}}}\sum_{n = 0}^{\infty} \frac{(q^{-n - m}, q)_m}{(q^{-1}; q^{-1})_m} x^{-n}
        \\ & = - \sum_{n = 0}^{\infty} \frac{(q^{-n-1}, q^{-1})_m}{(q;q)_m} x^{-n-m-1}
        \\ & = - \sum_{n < -m} \frac{(q^{n + m}, q^{-1})_m}{(q;q)_m} x^n
    \end{align*}
    The are two observations to make about this example. Firstly while the sum is $< -m$, it can be trivially extended to $< 0$ as the expression is exactly $0$ for $0 > n \geq -m$. Secondly, note that the initial equality is more commonly written as
    \[
        \sum_{n = 0}^{\infty} \frac{(q; q)_{m + n}}{(q;q)_n(q;q)_m} x^n = \frac{1}{(x;q)_{m + 1}}.
    \]
    Writing it in this way however, makes the ``inversion" difficult as as $\frac{1}{(q; q)_{-n}} = 0$ for all $n > 0$.
    
    The inverse state sum can be thought of as another example of this phenomenon. Given a knot $K$ with a chosen $n + 1$ strand braid diagram recall, as in Equation \ref{eq: r mat braid sum}, that we have a formal expansion:
    \[
        \sum_{\substack{\textbf{a}_1, \cdots, \textbf{a}_n \\ \textbf{c}_1, \cdots \textbf{c}_m}} \prod_{i = 1}^m x^{\frac{N - 1}{2}} q^{-|\textbf{a}_i|} \prod_{\alpha \in \text{crossings}} {{}_{\G{sl}_N}R^{e_{\alpha}}}_{\textbf{i}_\alpha, \textbf{j}_\alpha}^{\textbf{i}_\alpha', \textbf{j}_\alpha'}
    \]
    The idea is, as above, to invert some of these summations, choosing $2$ subsets $P \subset \{1, \cdots, n\}$ and $Q \subset \{1, \cdots, m\}$, and computing:
    \[
        (-1)^{|P| + |Q|} \sum_{\substack{-\textbf{a}_{i}, i\in P \\ -\textbf{c}_j, j \in Q}} \  \sum_{\substack{\textbf{a}_{i}, i\notin P \\ \textbf{c}_j, j \notin Q}} \prod_{i = 1}^m x^{\pm \frac{N - 1}{2}} q^{\pm |\textbf{a}_i|} \prod_{\alpha \in \text{crossings}} {{}_{\G{sl}_N}R^{e_{\alpha}}}_{\textbf{i}_\alpha, \textbf{j}_\alpha}^{\textbf{i}_\alpha', \textbf{j}_\alpha'}
    \]
    Here the $\pm$ in the $q$ power is $+$ if $i\in P$ and $-$ otherwise. Additionally, a negative state $-\textbf{a}$ is a sequence of integers $0 < a_n < \cdots < a_1$. Note that we use $<$ and not $\leq$. This is simply a reflection of our choice of basis and in the quantum polynomial representation this represents summing over monomials with strictly negative powers.
    
    In \cite{Park3} it was shown that, if it is possible to choose the inversions wisely such that this inverse state sum converges absolutely then the result is precisely the $F_K$ invariant that we are looking for.
    
    Computationally, what is going on is that when $\textbf{b}, \textbf{b}'$ are negative\footnote{The method in \cite{Park3} is slightly more general than this, but we ignore this here.}, we can still make sense of the ${{}_{\G{sl}_N}R^{-1}}_{\textbf{a}, \textbf{b}}^{\textbf{a}', \textbf{b}'}$ matrix elements and these will give a series in $x^{-1}$. Thus if we can give `inverted' \textbf{b} inputs to all $R^{-1}$ matrices, the state sum will converge. The surprising feature from the computational standpoint is that this inverted summation is exactly $F_K$.
    
    \subsection{The \texorpdfstring{$R^{-1}$}{R\^-1} matrix for inverted states}
        Let's study how to interpret the $R^{-1}$ matrix for negative states. For simplicity we mostly focus on $\G{sl}_{2,3}$ here but this easily extends to higher $N$. Starting with $\G{sl}_2$, each state $\ket{\textbf{a}}$ is defined by a single number $a \geq 0$ and so the $R^{-1}$ matrix elements coming from Equation \eqref{eq: RMatrix slN} specialise to:
        \[
            {{}_{\G{sl}_2}R^{-1}}_{\textbf{a}, \textbf{b}}^{\textbf{a}', \textbf{b}'} = (-1)^{b + a'}q^{\frac{b'}{2} + \frac{a^2}{2} - aa' - \frac{a}{2} - \frac{(b')^2}{2}} x^{\frac{1}{2}(b + b')}\frac{(q^{-b}; q)_{b - a'}(x^{-1} q^{a}; q)_{b - a'}}{(q^{-1};q^{-1})_{b - a'}}
        \]
        In comparison to the generic case, note that the summation over $r$ disappears\footnote{In general the ${}_{\G{sl}_N}R_{\textbf{a}_1, \textbf{b}_1}^{\textbf{a}_2, \textbf{b}_2}$ matrix element will be a sum over $\frac{(N - 1)(N - 2)}{2}$ $r$ variables.} as the only non $0$ term occurs at $r = a - b'$. Additionally, we have used the equality $(x q^b; q)_r = (-1)^rx^rq^{rb}q^{\frac{r(r-1)}{2}} (x^{-1} q^{-b}; q^{-1})_r$. Next, observe that
        \[
            \frac{(q^{-b}; q)_{b - a'}}{(q^{-1};q^{-1})_{b - a'}} = \frac{(q^{-1}, q^{-1})_b}{(q^{-1};q^{-1})_{b - a'}(q^{-1};q^{-1})_{a'}} = \frac{(q^{-b}; q)_{a'}}{(q^{-1};q^{-1})_{a'}}
        \]
        The beauty of this equality is that while all terms are equal for the usual $b > a' > 0$ situation, if $a' > 0 > b$ the right hand side still makes sense and similarly if $0 > b > a$ then we can use the left hand side. Mostly we will be concerned with the former situation and so use the $R^{-1}$ matrix
        \[
            (-1)^{b + a'}q^{\frac{b'}{2} + \frac{a^2}{2} - aa' - \frac{a}{2} - \frac{(b')^2}{2}} x^{\frac{1}{2}(b + b')}\frac{(q^{-b}; q)_{a'}(x^{-1} q^{a}; q)_{b - a'}}{(q^{-1};q^{-1})_{a'}}
        \]
        The point as mentioned earlier is that if $b, b'$ are less than $0$, this matrix will produce a series in $x^{-1}$ and this will allow us to produce convergent answers from the state sums in some cases.

        We take a similar approach for $\G{sl}_3$. In this case, each state $\ket{\textbf{a}}$ is defined by a pair of numbers $a_1 \geq a_2 \geq 0$ and so we get:
        \begin{align*}
            {{}_{\G{sl}_3}R^{-1}}_{\textbf{a}, \textbf{b}}^{\textbf{a}', \textbf{b}'} & = \sum_{r \geq 0} (-1)^{b_1 + a_1' + b_2 + a'_2 + r}q^{-C_3(r)}x^{\frac{1}{2}(b_1 + b'_1)}(x^{-1} q^{a_1}; q)_{b_1 - a'_1}
            \\ & \quad \quad \times \frac{(q^{-b_2};q)_{a'_2 + r}(q^{a_2 - a_1};q)_{b_2 - a'_2 - r}(q^{a'_2 - b_1 + r};q)_{a'_1 - a'_2}}{(q^{-1};q^{-1})_{a'_1 - a'_2}(q^{-1};q^{-1})_{a'_2}(q^{-1};q^{-1})_{r}}
            \\ C_3(r) & = \frac{1}{4}\Big(2(b_2^2 + {a'_2}^2 + r^2) + (4a_1 + 4r - a_2)a'_1 + (2 + a_1)a'_2 - 2(b - a')(b' + a' - 1)
            \\ & \quad \quad + b_1(a'_2 - a_2) + 2r - b_2(2 + 3a_1 - 4a_2 + a'_1 + 4a'_2 + 4r)\Big)
        \end{align*}
        Note that as mentioned previously, the summand is $0$ unless
        \[
            0, a_2 + b_2 - a_1 - a'_2 \leq r \leq b_1 - a'_1, b_2 - a'_2.
        \]
        and so once we specialise \textbf{a} and \textbf{b} this is a finite sum. Our work has essentially been done for us already here as these expressions already produce sensible answers for $0 < b_2 < b_1$.
        
        This method is very robust and continues to work for higher $\G{sl}_N$. It can also be combined with the stratified state sum method to produce predictions for far more knots. The only downside is that is is more computationally expensive and so we only present computations in the $N = 2, 3$ cases here.

\subsection{Homogeneous braid knots}
    A homogeneous braid is a braid such that for every $i$, either $\sigma_i$ or $\sigma_i^{-1}$ appears but never both. These are exactly the braids for which the inverse state sum technique  works most easily.
    
\begin{figure}[htp]
    	\centering
    	
    	\begin{tikzpicture}
    	
    	    \def\fsizeone{9}
    	    \def\fsizetwo{10}
    	    \def\size{0.75}
    	    \def\xshft{-2}
    	    
    	    
    	    \node at (-5 + \xshft, 0) {$4_1$ Knot: $\sigma_1\sigma_2^{-1}\sigma_1\sigma_2^{-1} \mapsto$};
    	    
    	    \node at (\xshft - \size/2, \size) {\fontsize{\fsizeone}{\fsizetwo}\selectfont \textbf{0}};
    	    \node at (\xshft - \size/2, 0) {\fontsize{\fsizeone}{\fsizetwo}\selectfont \textbf{a}};
    	    \node at (\xshft - \size/2, -\size) {\fontsize{\fsizeone}{\fsizetwo}\selectfont \textbf{-b}};
    	    
    	    \draw (\xshft + -\size/4, -\size) -- (\xshft + 5*\size/4, -\size);
    	    
    	    \HoroUnderCross[\xshft, 0](\size)
    	    
    	    \draw (\xshft + 5*\size/4, \size) -- (\xshft + 11*\size/4, \size);
    	    \node at (\xshft + 8*\size/4, 5*\size/4) {\fontsize{\fsizeone}{\fsizetwo}\selectfont \textbf{a}};
    	    
    	    \node at (\xshft + 5*\size/4, -\size/4) {\fontsize{\fsizeone}{\fsizetwo}\selectfont \textbf{0}};
    	    
    	    \HoroOverCross[\xshft + 3*\size/2, -\size](\size)
    	    
    	    \node at (\xshft + 11*\size/4, -\size/4) {\fontsize{\fsizeone}{\fsizetwo}\selectfont \textbf{0}};
    	    
    	    \draw (\xshft + 11*\size/4, -\size) -- (\xshft + 17*\size/4, -\size);
    	    \node at (\xshft + 14*\size/4, -5*\size/4) {\fontsize{\fsizeone}{\fsizetwo}\selectfont -\textbf{b}};
    	    
    	    \HoroUnderCross[\xshft + 6*\size/2, 0](\size)
    	    
    	    \node at (\xshft + 17*\size/4, \size/4) {\fontsize{\fsizeone}{\fsizetwo}\selectfont \textbf{a}};
    	    
    	    \draw (\xshft + 17*\size/4, \size) -- (\xshft + 23*\size/4, \size);
    	    
    	    \HoroOverCross[\xshft + 9*\size/2, -\size](\size)
    	    
    	    \node at (\xshft + 24*\size/4, \size) {\fontsize{\fsizeone}{\fsizetwo}\selectfont \textbf{0}};
    	    \node at (\xshft + 24*\size/4, 0) {\fontsize{\fsizeone}{\fsizetwo}\selectfont \textbf{a}};
    	    \node at (\xshft + 24*\size/4, -\size) {\fontsize{\fsizeone}{\fsizetwo}\selectfont \textbf{-b}};
    	    
    	    
    	    \def\yshft{-2.5}
    	    
    	    \node at (-5 + \xshft, \yshft) {$6_2$ Knot: $\sigma_1^3\sigma_2^{-1}\sigma_1\sigma_2^{-1} \mapsto$};
    	    
    	    \node at (\xshft - 4*\size/2, \size + \yshft) {\fontsize{\fsizeone}{\fsizetwo}\selectfont \textbf{0}};
    	    \node at (\xshft - 4*\size/2, 0 + \yshft) {\fontsize{\fsizeone}{\fsizetwo}\selectfont \textbf{a}};
    	    \node at (\xshft - 4*\size/2, -\size + \yshft) {\fontsize{\fsizeone}{\fsizetwo}\selectfont -\textbf{b}};
    	    
    	    \draw (\xshft + -7*\size/4, -\size + \yshft) -- (\xshft + 11*\size/4, -\size + \yshft);
    	    
    	    \HoroUnderCross[\xshft - 3*\size/2, 0 + \yshft](\size)
    	    
    	    \node at (\xshft - \size/4, 5*\size/4 + \yshft) {\fontsize{\fsizeone}{\fsizetwo}\selectfont \textbf{a}};
    	    \node at (\xshft - \size/4, -\size/4 + \yshft) {\fontsize{\fsizeone}{\fsizetwo}\selectfont \textbf{0}};
    	    
    	    \HoroUnderCross[\xshft, 0 + \yshft](\size)
    	    
    	    \node at (\xshft + 5*\size/4, 5*\size/4 + \yshft) {\fontsize{\fsizeone}{\fsizetwo}\selectfont \textbf{c}\textsubscript{1}};
    	    \node at (\xshft + 5*\size/4, -\size/4 + \yshft) {\fontsize{\fsizeone}{\fsizetwo}\selectfont \textbf{a}-\textbf{c}\textsubscript{1}};
    	    
    	    \HoroUnderCross[\xshft + 3*\size/2, 0 + \yshft](\size)
    	    
    	    \draw (\xshft + 11*\size/4, \size + \yshft) -- (\xshft + 17*\size/4, \size + \yshft);
    	    \node at (\xshft + 14*\size/4, 5*\size/4 + \yshft) {\fontsize{\fsizeone}{\fsizetwo}\selectfont \textbf{a}-\textbf{c}\textsubscript{2}};
    	    
    	    \node at (\xshft + 11*\size/4, -\size/4 + \yshft) {\fontsize{\fsizeone}{\fsizetwo}\selectfont \textbf{c}\textsubscript{2}};
    	    
    	    \HoroOverCross[\xshft + 6*\size/2, -\size + \yshft](\size)
    	    
    	    \node at (\xshft + 17*\size/4, -\size/4 + \yshft) {\fontsize{\fsizeone}{\fsizetwo}\selectfont \textbf{0}};
    	    
    	    \draw (\xshft + 17*\size/4, -\size + \yshft) -- (\xshft + 23*\size/4, -\size + \yshft);
    	    \node at (\xshft + 20*\size/4, -5*\size/4 + \yshft) {\fontsize{\fsizeone}{\fsizetwo}\selectfont \textbf{c}\textsubscript{2}-\textbf{b}};
    	    
    	    \HoroUnderCross[\xshft + 9*\size/2, 0 + \yshft](\size)
    	    
    	    \node at (\xshft + 23*\size/4, \size/4 + \yshft) {\fontsize{\fsizeone}{\fsizetwo}\selectfont \textbf{a}-\textbf{c}\textsubscript{2}};
    	    
    	    \draw (\xshft + 23*\size/4, \size + \yshft) -- (\xshft + 29*\size/4, \size + \yshft);
    	    
    	    \HoroOverCross[\xshft + 12*\size/2, -\size + \yshft](\size)
    	    
    	    \node at (\xshft + 30*\size/4, \size + \yshft) {\fontsize{\fsizeone}{\fsizetwo}\selectfont \textbf{0}};
    	    \node at (\xshft + 30*\size/4, 0 + \yshft) {\fontsize{\fsizeone}{\fsizetwo}\selectfont \textbf{a}};
    	    \node at (\xshft + 30*\size/4, -\size + \yshft) {\fontsize{\fsizeone}{\fsizetwo}\selectfont -\textbf{b}};
    	    
    	    
    	    \node at (-5 + \xshft, 2*\yshft) {$6_3$ Knot: $\sigma_1^2\sigma_2^{-1}\sigma_1\sigma_2^{-2} \mapsto$};
    	    
    	    \node at (\xshft - 4*\size/2, \size + 2*\yshft) {\fontsize{\fsizeone}{\fsizetwo}\selectfont \textbf{0}};
    	    \node at (\xshft - 4*\size/2, 0 + 2*\yshft) {\fontsize{\fsizeone}{\fsizetwo}\selectfont \textbf{a}};
    	    \node at (\xshft - 4*\size/2, -\size + 2*\yshft) {\fontsize{\fsizeone}{\fsizetwo}\selectfont -\textbf{b}};
    	    
    	    \draw (\xshft + -7*\size/4, -\size + 2*\yshft) -- (\xshft + 5*\size/4, -\size + 2*\yshft);
    	    
    	    \HoroUnderCross[\xshft - 3*\size/2, 0 + 2*\yshft](\size)
    	    
    	    \node at (\xshft - \size/4, 5*\size/4 + 2*\yshft) {\fontsize{\fsizeone}{\fsizetwo}\selectfont \textbf{a}};
    	    \node at (\xshft - \size/4, -\size/4 + 2*\yshft) {\fontsize{\fsizeone}{\fsizetwo}\selectfont \textbf{0}};
    	    
    	    \HoroUnderCross[\xshft, 0 + 2*\yshft](\size)
    	    
    	    \draw (\xshft + 5*\size/4, \size + 2*\yshft) -- (\xshft + 11*\size/4, \size + 2*\yshft);
    	    \node at (\xshft + 8*\size/4, 5*\size/4 + 2*\yshft) {\fontsize{\fsizeone}{\fsizetwo}\selectfont \textbf{a}-\textbf{c}\textsubscript{1}};
    	    
    	    \node at (\xshft + 5*\size/4, -\size/4 + 2*\yshft) {\fontsize{\fsizeone}{\fsizetwo}\selectfont \textbf{c}\textsubscript{1}};
    	    
    	    \HoroOverCross[\xshft + 3*\size/2, -\size + 2*\yshft](\size)
    	    
    	    \node at (\xshft + 11*\size/4, -\size/4 + 2*\yshft) {\fontsize{\fsizeone}{\fsizetwo}\selectfont \textbf{0}};
    	    
    	    \draw (\xshft + 11*\size/4, -\size + 2*\yshft) -- (\xshft + 17*\size/4, -\size + 2*\yshft);
    	    \node at (\xshft + 14*\size/4, -5*\size/4 + 2*\yshft) {\fontsize{\fsizeone}{\fsizetwo}\selectfont \textbf{c}\textsubscript{1}-\textbf{b}};
    	    
    	    \HoroUnderCross[\xshft + 6*\size/2, 0 + 2*\yshft](\size)
    	    
    	    \node at (\xshft + 17*\size/4, -\size/4 + 2*\yshft) {\fontsize{\fsizeone}{\fsizetwo}\selectfont \textbf{a}-\textbf{c}\textsubscript{1}};
    	    
    	    \node at (\xshft + 23*\size/4, -5*\size/4 + 2*\yshft) {\fontsize{\fsizeone}{\fsizetwo}\selectfont \textbf{c}\textsubscript{2}-\textbf{b}};
    	    
    	    \HoroOverCross[\xshft + 9*\size/2, -\size + 2*\yshft](\size)
    	    
    	    \node at (\xshft + 24*\size/4, \size/4 + 2*\yshft) {\fontsize{\fsizeone}{\fsizetwo}\selectfont \textbf{a}-\textbf{c}\textsubscript{2}};
    	    
    	    \draw (\xshft + 17*\size/4, \size + 2*\yshft) -- (\xshft + 29*\size/4, \size + 2*\yshft);
    	    
    	    \HoroOverCross[\xshft + 12*\size/2, -\size + 2*\yshft](\size)
    	    
    	    \node at (\xshft + 30*\size/4, \size + 2*\yshft) {\fontsize{\fsizeone}{\fsizetwo}\selectfont \textbf{0}};
    	    \node at (\xshft + 30*\size/4, 0 + 2*\yshft) {\fontsize{\fsizeone}{\fsizetwo}\selectfont \textbf{a}};
    	    \node at (\xshft + 30*\size/4, -\size + 2*\yshft) {\fontsize{\fsizeone}{\fsizetwo}\selectfont -\textbf{b}};

        \end{tikzpicture}
        \captionsetup{singlelinecheck=off}
    	\caption[.]{Labeled Braid diagrams for the $4_1$ and $6_2$ and $6_3$ knots.}
    	\label{fig: Braid Diagrams 41 62 63}
    \end{figure}
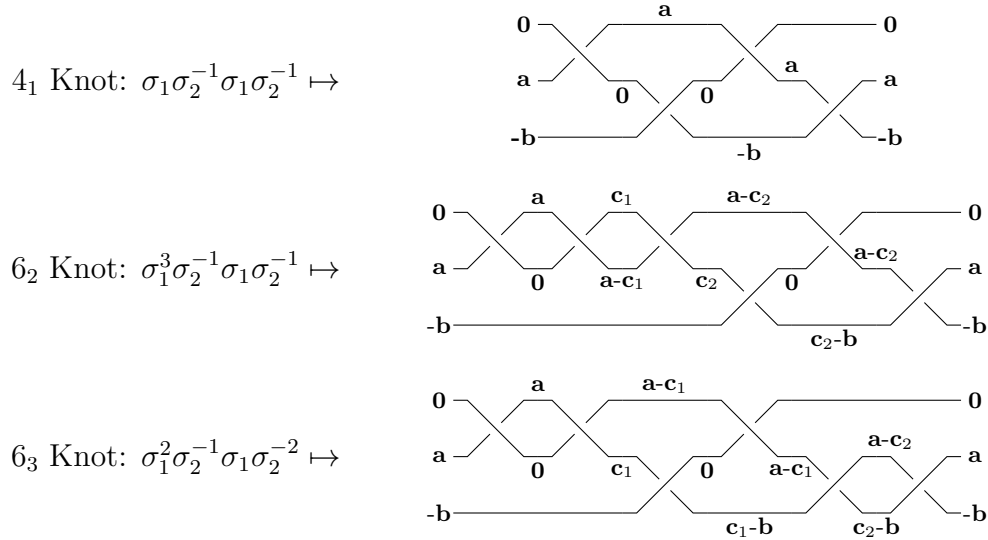

    The simplest knots which have homogeneous braid representatives are the $4_1, 6_2$ and $6_3$ knots. Labeled braid diagrams for these knots are in Figure \ref{fig: Braid Diagrams 41 62 63} and the corresponding $F^K_N$ are given in Tables \ref{tab: FK 41}, \ref{tab: FK 62} and \ref{tab: FK 63}. The $N = 2$ results previously appeared in \cite{Park3} but we give them here both for a consistency check\footnote{As usual, to compare these results we first need to align conventions.} and to give hints towards possible $a$ deformations.
    
    It can be easily checked that the expressions in Table \ref{tab: FK 41} match the $a = q^2, q^3$ specialisations of the $F_{4_1}(x, a, q)$ given in \cite{EGGKPSS}. Similarly the $\G{sl}_2$ rows in Tables \ref{tab: FK 62}, \ref{tab: FK 63} match the computations in \cite{Park3}. The $6_2, 6_3$ examples are also interesting as they are the first examples where the $a$ deformation is still unknown. Thus the $\G{sl}_3$ examples here should serve as an important cross check for future work in that direction. Unfortunately, there is not enough data here to provide a sensible guess for what the $a$ deformations might be.
    
    \begin{table}
        \centering
        \renewcommand{\arraystretch}{1.5}
        \begin{tabular}{ c | m{15cm} }
        & \multicolumn{1}{c}{$F^+_{4_1}(x, q)$} \\ \hline
        \shortstack{$\G{sl}_2$} & $1 + 3qx + q(1 + 6q + q^2) x^2 + q(2 + 3q + 11q^2 + 3q^3 + 2q^4) x^3 + (1 + 3q + 6q^2 + 8q^3 + 19q^4 + 8q^5 + 6q^6 + 3q^7 + q^8)x^4 + (2q^{-1} + 2 + 7q + 10q^2 + 16q^3 + 18q^4 + 34q^5 + 18q^6 + 16q^7 + 10q^8 + 7q^9 + 2q^{10} + 2q^{11})x^5$ \\  \hline
        \shortstack{$\G{sl}_3$} & $1 + 3q (1 + q) x + q(1 + 7q + 9q^2 + 7q^3 +q^4) x^2 + q (2 + 5q + 16q^2 + 22q^3 + 22q^4 + 16q^5 + 5q^6 + 2q^7) x^3 + (1 + 4q + 10q^2 + 18q^3 + 38q^4 + 51q^5 + 56q^6 + 51q^7 + 38q^8 + 18q^9 + 10q^{10} + 4q^{11} + q^{12}) x^4 + (2 q^{-1} + 4 + 11q + 21q^2 + 37q^3 + 57q^4 + 92q^5 + 119q^6 + 134q^7 + 134q^8 + 119q^9 + 92q^{10} + 57q^{11} + 37q^{12} + 21q^{13} + 11q^{14} + 4q^{15} + 2q^{16}) x^5$ \\ \hline
        \end{tabular}
        \caption{$F_K$ invariant for the $4_{1}$ knot. This can be matched against the specialisations of the $a$ deformed series given in \cite{EGGKPSS}, Equation (38)}
        \label{tab: FK 41}
    \end{table}
    
    \begin{table}
        \centering
        \renewcommand{\arraystretch}{1.5}
        \begin{tabular}{ c | m{15cm} }
        & \multicolumn{1}{c}{$F^+_{6_2}(x, q)$} \\ \hline
        \shortstack{$\G{sl}_2$} & $1 + 3qx + (q + 6q^2 - q^3)x^2 + q(2 + 3q + 10q^2 - 3q^3)x^3 + (1 + 3q + 6q^2 + 6q^3 + 14q^4 - 6q^5 + 2q^7)x^4 + (2q^{-1} + 2 + 7q + 9q^2 + 12q^3 + 8q^4 + 18q^5 - 10q^6 + q^7 + 6q^8 + 2q^9)x^5$ \\  \hline
        \shortstack{$\G{sl}_3$} & $1 + 2q(1 + 2q)x + q(2 + 3q + 8q^2 + 9q^3 -q^4)x^2 + 2q(2 + 2q + 6q^2 + 6q^3 + 9q^4 + 7q^5 - 2q^6)x^3 + (3 + 6q + 8q^2 + 22q^3 + 20q^4 + 35q^5 + 25q^6 + 29q^7 + 17q^8 - 9q^9 +q^{10} + 3q^{11})x^4 + (6q^{-1} + 6 + 14q + 24q^2 + 36q^3 + 38q^4 + 66q^5 + 50q^6 + 67q^7 + 37q^8 + 38q^9 + 19q^{10} - 13q^{11} + 10q^{12} + 11q^{13} + 5q^{14})x^5$ \\  \hline
        \end{tabular}
        \caption{$F_K$ invariant for the $4_{1}$ knot and symmetric series of representations on $\G{sl}_N$ for small $N$.}
        \label{tab: FK 62}
    \end{table}
    
    \begin{table}
        \centering
        \renewcommand{\arraystretch}{1.5}
        \begin{tabular}{ c | m{15cm} }
        & \multicolumn{1}{c}{$F^+_{6_3}(x, q)$} \\ \hline
        \shortstack{$\G{sl}_2$} & $1 + 3qx - q(1 - 6q + q^2)x^2 -(-1 +q)^2q (2 + 7q + 2q^2)x^3 - 2q (1 + 3q + 3q^2 - 8q^3 + 3q^4 + 3q^5 +q^6)x^4 + (2 - 2q - 6q^2 - 12q^3 - 8q^4 + 25q^5 - 8q^6 - 12q^7 - 6q^8 - 2q^9 + 2q^{10})x^5$ \\  \hline
        \shortstack{$\G{sl}_3$} & $1 + 3q (1 + q)x - q(1 - 5q - 9q^2 - 5q^3 +q^4)x^2 - q (2 + 5q - 4q^2 - 15q^3 - 15q^4 - 4q^5 + 5q^6 + 2q^7)x^3 - q (2 + 8q + 17q^2 + 6q^3 - 16q^4 - 26q^5 - 16q^6 + 6q^7 + 17q^8 + 8q^9 + 2q^{10})x^4 + (2 - 6q^2 - 22q^3 - 43q^4 - 32q^5 + 4q^6 + 34q^7 + 34q^8 + 4q^9 - 32q^{10} - 43q^{11} - 22q^{12} - 6q^{13} + 2q^{15})x^5$ \\ \hline
        \end{tabular}
        \caption{$F_K$ invariant for the $6_{3}$ knot and symmetric series of representations on $\G{sl}_N$ for small $N$.}
        \label{tab: FK 63}
    \end{table}

\section{Future Directions}

The work presented here naturally leaves some questions to be answered.

\begin{enumerate}
    \item Can the proof of the inverse state sum method be extended to work in the $\G{sl}_N$ case? 
    \item Is there a general construction which will work to compute $F^N_K$ for any integer $N$ and knot $K$?  
    \item Can this construction be generalised to other infinite families of representations?
    \item Can we find a multi-variable large colour $R$ matrix for $\G{sl}_3$. It should have $3$ variables $x_1, x_2, q$ and there should be a specialisation which recovers the symmetrically coloured $R$ matrix given here?
    \item Is there a way to pass pass from an $R$ matrix summation to a quiver form which does not grow as $N$ increases?
    \item Is there a state sum expression which can generate the $a$ deformed $F_K$ invariant $F_K(x, a, q)$?
\end{enumerate}

\section{Acknowledgements}
    Many thanks to Sergei Gukov, Sunghyuk Park and Piotr Kucharski for useful conversations. I was supported by the National Science Foundation under Grant No. NSF DMS 1664227.

\bibliography{MCSV}

\newcommand{\etalchar}[1]{$^{#1}$}
\begin{thebibliography}{EGG{\etalchar{+}}22b}

\bibitem[BNG]{BNG}
Dror Bar-Natan and Stavros Garoufalidis.
\newblock On the {Melvin-Morton-Rozansky} conjecture.
\newblock {\em Invent. Math.}, 125(1):103--133, 1996.





\bibitem[Bur]{Bur}
Nigel Burroughs.
\newblock {The} universal {R-matrix} for {Uqsl(3)} and beyond!
\newblock {\em Communications in Mathematical Physics}, 127(1):109--128, 1990.





\bibitem[EGGKPS]{EGGKPS}
Tobias Ekholm, Angus Gruen, Sergei Gukov, Piotr Kucharski, Sunghyuk Park, and
  Piotr Su\l{}kowski.
\newblock $\widehat{Z}$ at large $n$: from curve counts to quantum modularity.
\newblock {\em Communications in Mathematical Physics}, 396, 08 2022.
\newblock arXiv:2005.13349.





\bibitem[EGGKPSS]{EGGKPSS}
Tobias Ekholm, Angus Gruen, Sergei Gukov, Piotr Kucharski, Sunghyuk Park, Marko
  Stošić, and Piotr Sułkowski.
\newblock Branches, quivers, and ideals for knot complements.
\newblock {\em Journal of Geometry and Physics}, 177:104520, 2022.





\bibitem[FGS]{FGS}
Hiroyuki Fuji, Sergei Gukov, and Piotr Sulkowski.
\newblock Super-{$A$}-polynomial for knots and {BPS} states.
\newblock {\em Nucl. Phys. B}, 867:506, 2013.
\newblock arXiv:1205.1515.





\bibitem[GM]{GM}
Sergei Gukov and Ciprian Manolescu.
\newblock A two-variable series for knot complements, 2019.
\newblock arXiv:1904.06057.





\bibitem[Gar]{Gar}
Stavros Garoufalidis.
\newblock On the charactersitic and deformation varieties of a knot.
\newblock {\em Geometry and Topology Monographs}, 7:291--304, 2004.
\newblock math/0306230.





\bibitem[Guk]{Guk}
Sergei Gukov.
\newblock Three-dimensional quantum gravity, {Chern-Simons} theory, and the
  {A-polynomial}.
\newblock {\em Commun. Math. Phys.}, 255:577--627, 2005.
\newblock hep-th/0306165.





\bibitem[Jon]{Jon}
Vaughan Jones.
\newblock A polynomial invariant for knots via von {Neumann} algebras.
\newblock {\em Bull. Amer. Math. Soc.}, 12:103--111, 1985.





\bibitem[KRSS1]{KRSS1}
Piotr Kucharski, Markus Reineke, Marko Stosic, and Piotr Sulkowski.
\newblock {BPS} states, knots and quivers.
\newblock {\em Phys. Rev.}, D96(12):121902, 2017.
\newblock arXiv:1707.02991.





\bibitem[KRSS2]{KRSS2}
Piotr Kucharski, Markus Reineke, Marko Stosic, and Piotr Sulkowski.
\newblock Knots-quivers correspondence.
\newblock {\em Adv. Theor. Math. Phys.}, 23(7):1849--1902, 2019.
\newblock arXiv:1707.04017.





\bibitem[Kuch]{Kuch}
Piotr Kucharski.
\newblock {Quivers for 3-manifolds: the correspondence, BPS states, and 3d $
  \mathcal{N} $ = 2 theories}.
\newblock {\em JHEP}, 09:075, 2020.
\newblock arXiv:2005.13394.





\bibitem[LM]{LM}
Charles Livingston and Allison~H. Moore.
\newblock Knotinfo: Table of knot invariants.
\newblock URL: \url{knotinfo.math.indiana.edu}, 3 2022.





\bibitem[MM]{MM}
P.~M. Melvin and H.~R. Morton.
\newblock The coloured {Jones} function.
\newblock {\em Comm. Math. Phys.}, 169(3):501--520, 1995.





\bibitem[Park1]{Park1}
Sunghyuk Park.
\newblock {Higher rank $\hat{Z}$ and $F_K$}.
\newblock {\em {SIGMA}}, {16}(044), 2020.
\newblock arXiv:1909.13002.





\bibitem[Park2]{Park2}
Sunghyuk Park.
\newblock {Large color R-matrix for knot complements and strange identities}.
\newblock {\em Journal of Knot Theory and Its Ramifications}, 29(14):2050097,
  Dec 2020.
\newblock arXiv:2004.02087.





\bibitem[Park3]{Park3}
Sunghyuk Park.
\newblock Inverted state sums, inverted habiro series, and indefinite theta
  functions.
\newblock {\em to appear}, 2021.





\bibitem[Roz]{Roz}
Lev Rozansky.
\newblock Higher order terms in the {Melvin-Morton} expansion of the colored
  {Jones} polynomial.
\newblock {\em Comm. Math. Phys.}, 183(2):291--306, 1997.





\bibitem[Tur]{Tur}
V.~G. Turaev.
\newblock The yang-baxter equation and invariants of links.
\newblock {\em Inventiones mathematicae}, 92:527 -- 553, 1988.





\bibitem[Wit]{Wit}
Edward Witten.
\newblock Quantum field theory and the {Jones} polynomial.
\newblock {\em Comm. Math. Phys.}, 121(3):351--399, 1989.





\end{thebibliography}
\bibliographystyle{abstract}

\end{document}